\input  amstex
\input amsppt.sty
\magnification1200
\vsize=23.5truecm
\hsize=16.5truecm
%\vcorrection{-10truemm}
\NoBlackBoxes

\def\supp{\operatorname{supp}}

\def\rnp{{\Bbb R}^n_+}

\def\rnpm{\Bbb R^n_\pm}
\def\crnp{\overline{\Bbb R}^n_+}

\def\crnpm{\overline{\Bbb R}^n_\pm}
\def\comega{\overline\Omega }

\def\Rn{\Bbb R^n}
\def\ang#1{\langle {#1} \rangle}

\def\Pfrac{\tsize\frac1{\raise 1pt\hbox{$\scriptstyle p$}}}
\def\pfrac{\frac1{\raise 1pt\hbox{$\scriptscriptstyle p$}}}
\def\Pfracc#1{\tsize\frac{#1}{\raise 1pt\hbox{$\scriptstyle p$}}}
\def\pfracc#1{\frac{#1}{\raise 1pt\hbox{$\scriptscriptstyle p$}}}

\def\simto{\overset\sim\to\rightarrow}

\def\Zfrac{\tsize\frac1{\raise 1pt\hbox{$\scriptstyle z$}}}
\def\zfrac{\frac1{\raise 1pt\hbox{$\scriptscriptstyle z$}}}

\def\rp{ \Bbb R_+}

\def\OP{\operatorname{OP}}

\def\R{\Bbb R}

\def\ol{\overline}
\def\SD{\Cal S}
\def\E{\Cal E}
\def\F{\Cal F}
\def\D{\Cal D}

\document
\topmatter
\title
 Regularity in $L_p$ Sobolev spaces of solutions to fractional heat equations 
\endtitle
\author Gerd Grubb \endauthor
\affil
{Department of Mathematical Sciences, Copenhagen University,
Universitetsparken 5, DK-2100 Copenhagen, Denmark.
E-mail {\tt grubb\@math.ku.dk}}\endaffil
\rightheadtext{Fractional heat equations}
\abstract
This work contributes in two areas, with sharp results, to the current investigation of
regularity of solutions of heat equations with a nonlocal operator $P$:
$$
Pu+\partial_tu=f(x,t), \text{ for  }x\in\Omega \subset{\Bbb R}^n,\;t\in
I\subset {\Bbb R}.\tag*
$$ 

1) For strongly elliptic
pseudodifferential operators ($\psi $do's) $P$ on ${\Bbb R}^n$ of
order $d\in\rp$, a symbol calculus on ${\Bbb R}^{n+1}$ is introduced
that allows showing optimal regularity results, globally over ${\Bbb
R}^{n+1}$ and locally over $\Omega \times I$:
$$
f\in  H^{(s,s/d)}_{p,\operatorname{loc}}(\Omega \times I)\implies u\in  H^{(s+d,s/d+1)}_{p,\operatorname{loc}}(\Omega \times I),
$$
for
$s\in{\Bbb R}$, $1<p<\infty $. The $H_p^{(s,s/d)}$ are
anisotropic Sobolev spaces of Bessel-potential type, and there is a
similar result for Besov spaces.

2) Let $\Omega $ be smooth bounded, and let $P$ equal $(-\Delta )^a$
($0<a<1$), or its generalizations to singular integral operators with
regular kernels, generating stable L\'evy processes. 
With the Dirichlet condition
$\supp u\subset \comega$,  the
initial condition $u|_{t=0}=0$, and $f\in L_p(\Omega \times I)$, 
(*) has a unique solution 
$u\in L_p(I, H_p^{a(2a)}(\comega))$ with $\partial_tu\in L_p(\Omega
\times I)$. Here $H_p^{a(2a)}(\comega)=\dot H_p^{2a}(\comega)$
 if $a<1/p$, and is contained in $\dot H_p^{2a-\varepsilon }(\comega)$ if $a=1/p$,
but contains nontrivial elements from
$ d^a\ol H_p^{a}(\Omega )$ if $a>1/p$ (where $d(x)=
\operatorname{dist}(x,\partial\Omega )$). The interior regularity of $u$ is lifted when $f$ is more smooth.

\endabstract
\keywords Nonlocal heat equation; fractional Laplacian; stable
process; pseudodifferential operator; Dirichlet condition; Lp Sobolev regularity \endkeywords

\subjclass 35K05, 35K25, 47G30, 60G52 \endsubjclass
\endtopmatter

\head 0. Introduction \endhead

There is currently a great interest for evolution problems (heat equations)
$$
Pu(x,t)+\partial_tu(x,t)=f(x,t) \text{ on }\Omega \times I,\;\Omega
\text{ open }
\subset \R^n,\; I=\,]0,T[\,,\tag0.1
$$
where $P$ is a nonlocal operator, as for example the fractional
Laplacian $(-\Delta )^a$ ($0<a<1$) or other pseudodifferential
operators ($\psi $do's) or singular integral operators.
For
differential operators $P$, there are classical treatises such as
Ladyzhenskaya, Solonnikov and Uraltseva \cite{LSU68} with
$L_p$-methods, Lions and Magenes \cite{LM68} with $L_2$-methods,
Friedman \cite{F69} with $L_2$ semigroup methods, and numerous more
recent studies.
Motivated by the linearized Navier-Stokes problem, which can be
reduced to the form (0.1) with nonlocal ingredients,
the author jointly with Solonnikov treated such problems in
\cite{GS90}
(for $L_2$-spaces) and \cite{G95} (for $L_p$-spaces). In those papers, the operator $P$ fits
into the 
Boutet de Monvel calculus \cite{B71,G86,G90,G09}, and is necessarily of
integer order.

This does not cover fractional order operators, and 
the present paper aims to find techniques to handle (0.1) in
fractional cases. Firstly, we treat $\psi $do's without
boundary conditions in Sections 2 and 3, where we introduce a
systematic calculus that allows
showing regularity results globally in ${\Bbb R}^{n+1}$,  and
locally in arbitrary open subset $\Sigma \subset{\Bbb R}^{n+1}$, 
in terms of anisotropic function spaces
 described in detail in the Appendix:

\proclaim{Theorem 0.1} Let $P$ be a classical strongly elliptic $\psi $do  $P=\OP(p(x,\xi
))$ on ${\Bbb R}^n$ of order $d\in\rp$. Let $s\in{\Bbb R}$, $1<p<\infty
$. Then $P+\partial_t$ maps $H^{(s+d,s/d+1)}_p(\Rn\!\times\!\Bbb R)\to
H^{(s,s/d)}_p(\Rn\!\times\!\Bbb R)$.

$1^\circ$ Let $u\in
H^{(r,r/d)}_p(\Rn\!\times\!\Bbb R) $ for some large negative $r$ (this
holds in
particular if $u\in\E'({\Bbb R}^{n+1})$ or e.g.\ $L_p({\Bbb
R},\E'({\Bbb R}^n))$). Then 
$$
(P+\partial_t)u\in H^{(s,s/d)}_p(\Rn\!\times\!\Bbb
R)\implies
u\in H^{(s+d,s/d+1)}_p(\Rn\!\times\!\Bbb R).\tag0.2
$$

$2^\circ$ Let $\Sigma $ be an open subset of ${\Bbb
R}^{n+1}$, and let $u\in
H^{(s,s/d)}_p(\Rn\!\times\!\Bbb R) $. Then
$$
(P+\partial_t)u|_{\Sigma }\in H^{(s,s/d)}_{p,\operatorname{loc}}(\Sigma )\implies
u\in H^{(s+d,s/d+1)}_{p,\operatorname{loc}}(\Sigma ).
\tag0.3
$$ 
\endproclaim

The analogous result holds in Besov-spaces $B^{(s,s/d)}_p$, and there
is also a result in anisotropic H\"older spaces that can be derived
from (0.3) by letting $p\to\infty $.

A celebrated example to which the above theorem applies is
the fractional Laplacian 
$$
\aligned
(-\Delta )^au&=\operatorname{Op}(|\xi |^{2a})u=
\Cal F^{-1}(|\xi |^{2a}\hat u(\xi )),\text{ also defined by}\\
(-\Delta )^au(x)&=c_{n,a}PV\int_{{\Bbb
R}^n}\frac{u(x)-u(y)}{|x-y|^{n+2a}}
\,dy,
\endaligned\tag0.4
$$
which has interesting applications in
probability, finance, mathematical physics and differential geometry.
Along with (0.4), one considers more general translation-invariant singular integral operators
$$
Pu(x)=PV\int_{{\Bbb
R}^n}\frac{(u(x)-u(x+y))k(y/|y|)}{|y|^{n+2a}}
\,dy,\tag0.5
$$
with $k(y)$ even and positive on $S^{n-1}$ (cf.\ e.g.\ the survey of
Ros-Oton \cite{R16}); they are infinitesimal generators of stable
L\'evy processes. Further generalizations with nonhomogeneous or
nonsmooth kernels are also studied. When $k\in C^\infty $,  
the operator (0.5) is a $\psi $do of order $2a$ with positive homogeneous symbol
$p(\xi )$ even in $\xi $, and
Theorem 0.1  applies with $d=2a$.

We underline that the above regularity results apply not only to such operators, but also to $x$-dependent
$\psi $do's,
and to $\psi $do's with complex symbol,  without special symmetries and with a different 
behavior at boundaries (no ``transmission property'').
(An example is the square-root
Laplacian with drift $(-\Delta )^{\frac12}+b(x)\cdot \nabla$.) 

For (0.5), regularity questions for solutions of (0.1) have been treated recently by Leonori, Peral, Primo and
Soria \cite{LPPS15} in $L_r(I,L_q(\Omega ))$-spaces, by Fernandez-Real
and Ros-Oton \cite{FR17} in anisotropic H\"older spaces, and by
Biccari, Warma and Zuazua \cite{BWZ17} for $(-\Delta )^a$ in $L_p$-spaces valued in
local Sobolev spaces over
$\Omega $. Earlier results are shown e.g.\ in Felsinger and Kassmann
\cite{FK13} and Chang-Lara and Davila
\cite{CD14} (H\"older properties), and Jin and Xiong \cite{JX15}
(Schauder estimates); the references in the mentioned works give
further information, also on related heat kernel estimates.

The second aim of our paper is to obtain a global result in $L_p$
Sobolev spaces for the
heat equation (0.1) for 
$(-\Delta  )^a$ or (0.5), with Dirichlet boundary condition on a
bounded smooth domain
$\Omega $, giving a detailed description of the solution.
By combining the characterization of the Dirichlet domain obtained in \cite{G15} with a semigroup
theorem of Lamberton \cite{L87} put forward in \cite{BWZ17}, we show
in Section 4:

\proclaim{Theorem 0.2} Let $1<p<\infty $. When $P=(-\Delta )^a$, or is an operator as in
{\rm (0.5)} with $k\in C^\infty $, and  $\Omega \subset {\Bbb R}^n$ is bounded
smooth, then the  evolution problem 
$$
\aligned
Pu+\partial_tu&=f\text{ on }\Omega \times I,\\
u&=0\text{ on }(\R^n\setminus\Omega) \times  I,\\
u&=0\text{ for }t=0,
\endaligned \tag 0.6
$$
has for any $f\in L_p(\Omega\times  I)$ a unique solution $u(x,t)\in
C^0(\overline I,L_p(\Omega ))$, that satisfies:
$$
u\in L_p(I,H_p^{a(2a)}(\comega))\cap H^1_p( I,L_p(\Omega )).\tag0.7
$$
Here $H_p^{a(2a)}(\comega)$ is the domain of the Dirichlet $L_p$-realization 
$P_{\operatorname{Dir},p}$ of $P$ and equals\linebreak $\dot H_p^{2a}(\comega)+V$,
where $V=0$ if $a<1/p$, $V\subset \dot H_p^{2a-\varepsilon }(\comega)$
if 
$a=1/p$, $V\subset d^a\ol H_p^{a}(\Omega )$ if $a>1/p$ (here $d(x)= \operatorname{dist}(x,\partial\Omega )$).
\endproclaim

For the time-dependent
 problem, this precision is new. 

An application of Theorem 0.1 $2^\circ$ gives moreover:

\proclaim{Corollary 0.3} Let $u$ be as in Theorem {\rm 0.2}, and let
$r=2a$ if $a<1/p$, $r= a+1/p-\varepsilon $ if $a\ge
1/p$ (for some small $\varepsilon >0$). Then for $0<s\le r$, 
$$f\in  H_{p,\operatorname{loc}}^{(s,s/(2a))}(\Omega \times I) \implies  
u\in  H_{p,\operatorname{loc}}^{(s+2a,s/(2a)+1)}(\Omega \times I).\tag0.8$$  
\endproclaim

For larger $s$, the local regularity (0.8) can be obtained via 
Theorem 0.1 $2^\circ$ if
one knows on beforehand that  $u\in  H_{p}^{(s,s/(2a))}({\Bbb R}^n\times I)$.
\medskip

\noindent {\it Plan of the paper.} Section 1 gives some definitions
and prerequisites.  The definitions and properties of anisotropic Sobolev spaces (of
Bessel-potential and Besov type) are collected in the Appendix. In Section 2 we show how an anisotropic  symbol 
calculus can be introduced, that covers the operators $P+\partial_t$
with a pseudodifferential $P$ of order $d\in \rp$. Section 3 gives the
proofs of the global and local regularity stated in Theorem 0.1.
In Section 4 we start by introducing some further prerequisites needed for
the global results on a bounded smooth subset, and then give the proof
of Theorem 0.2.

\head 1. Preliminaries \endhead

We shall use the notation set up in \cite{G15}, also used in
\cite{G14,G17}, and will just list
some important points here.

The function $\ang\xi $ equals $(|\xi |^2+1)^{\frac12}$. The
Fourier transform $\F$ is defined by $\hat
u(\xi )=\Cal F
u(\xi )= \int_{{\Bbb R}^n}e^{-ix\cdot \xi }u(x)\, dx$; it maps the
Schwartz space $\Cal S({\Bbb R}^n)$ of rapidly decreasing $C^\infty
$-functions into itself, and extends by duality to the temperate distributions  $\Cal S'({\Bbb R}^n)$. 

A {\it pseudodifferential operator} ($\psi $do) $P$ on ${\Bbb R}^n$ is
defined from a symbol $p(x,\xi )$ on ${\Bbb R}^n\times{\Bbb R}^n$ by 
$$
Pu=\operatorname{OP}(p(x,\xi ))u 
=(2\pi )^{-n}\int e^{ix\cdot\xi
}p(x,\xi )\hat u\, d\xi =\Cal F^{-1}_{\xi \to x}(p(x,\xi )\F u(\xi
)),\tag 1.1
$$  
using the Fourier transform $\F$. 
We refer to
textbooks such as H\"ormander \cite{H85}, Taylor \cite{T81}, Grubb \cite{G09} for the rules of
calculus (in particular the definition by oscillatory integrals in \cite{H85}). The symbols
$p$ of order $m\in{\Bbb R}$ we shall use  are generally taken to lie in the symbol space $S^m_{1,0}({\Bbb R}^n\times{\Bbb R}^n)$, consisting of
$C^\infty $-functions $p(x,\xi )$
such that $\partial_x^\beta \partial_\xi ^\alpha p(x,\xi
)$ is $O(\ang\xi ^{m-|\alpha |})$ for all $\alpha ,\beta $, for some
$m\in{\Bbb R}$, with global estimates for $x\in{\Bbb R}^n$ (as in
\cite{H85} start of Sect.\ 18.1, and \cite{G96}).  $P$ (of order $m$)
then  maps $H^s_p({\Bbb R}^n)$ continuously into
$H^{s-m}_p ({\Bbb R}^n)$ for all $s\in{\Bbb R}$ (cf.\ (1.4)).
$P$ is said to be {\it classical} when $p$ moreover 
has an asymptotic expansion $p(x,\xi )\sim \sum_{j\in{\Bbb
N}_0}p_j(x,\xi )$ with $p_j$ homogeneous in $\xi $ of degree $m-j$ for
$|\xi |\ge 1$, all $j$, and $$
p(x,\xi )- {\sum}_{j<J}p_j(x,\xi )\in S^{m-J}_{1,0}({\Bbb
R}^n\times \R^n),\text{ for all }J.\tag1.2
$$
 $P$ is then said to be (uniformly)  {\it elliptic}
when $p_0(x,\xi )\ge c|\xi |^m $ for $|\xi |\ge 1$, with $c>0$. To these operators one can add
the smoothing operators (mapping any  $H^s_p({\Bbb R}^n)$ into
$\bigcap_tH^t_p ({\Bbb R}^n)$), regarded as operators of order $-\infty $.
$S^m_{1,0}({\Bbb R}^n\times{\Bbb R}^n)$ will also be written
$S^m_{1,0}({\Bbb R}^{2n})$ for short.

Formula (1.1) will also be used in some cases of more general 
functions $p$ for which the formula can be given a sense, for example in case of the symbol $(\ang{\xi '}+i\xi
_n)^t$ in (4.2), the anisotropic symbols in 
Definition 2.1, the symbol $|\xi |^a$ in Example 2.10.

Recall  the composition rule: When $PQ=R$, then $R$ has
a symbol $r(x,\xi )$ with the following asymptotic expansion, called the
Leibniz product:
$$
r(x,\xi )\sim p(x,\xi )\# q(x,\xi )= \sum_{\alpha \in{\Bbb
N}_0^n}D_\xi ^\alpha p(x,\xi ) \partial_x^aq(x,\xi )/\alpha !.\tag 1.3
$$

We shall also define $\psi $do's on $\R^{n+1}$ with
variables denoted $(x,t)$, the dual variables denoted $(\xi ,\tau
)$. The symbols $h(x,t,\xi ,\tau )$ may satisfy other types of estimates
with repect to $(\xi ,\tau )$ than the $S^m_{1,0}$ estimates mentioned
above. To distinguish between operators on $\R^n$ and $\R^{n+1}$, we
may write  $\OP_x(p(x,\xi )
)$ resp.\ $\OP_{x,t}(h(x,t,\xi ,\tau ))$.

The standard Sobolev spaces $W^{s,p}({\Bbb R}^n)$, $1<p<\infty $ and
$s\ge 0$, have a different character according to whether $s$ is
integer or not. Namely, for $s$ integer, they consist of
$L_p$-functions with derivatives in $L_p$ up to order $s$, hence
coincide with the Bessel-potential spaces $H^s_p({\Bbb R}^n)$, defined
for $s\in{\Bbb R}$ by 
$$
H_p^s(\R^n)=\{u\in \SD'({\Bbb R}^n)\mid \F^{-1}(\ang{\xi }^s\hat u)\in
L_p(\R^n)\}.\tag1.4
$$
For noninteger $s$, the $W^{s,p}$-spaces coincide with the Besov spaces, defined e.g.\ 
as follows: For $0<s<2$,
$$
f\in
B^s_p({\Bbb R}^{n})\iff \|f\|_{L_p}^p+ \int_{\R^{2n}}\frac{|f(x)+f(y)-2f((x+y)/2)|^p}{|x+y|^{n+ps}}\,dxdy<\infty ;\tag1.5
$$
and $B^{s+t}_p(\R^n)=(1-\Delta )^{-t/2}B^s_p({\Bbb R}^n)$ for all $t\in{\Bbb R}$.
The Bessel-potential spaces are important because they are most
directly related to $L_p(\R^n)$; the Besov spaces have other
convenient properties, and are
needed for boundary value problems in an $H^s_p$-context,
 because they are the correct range spaces for
trace maps (both from $H^s_p$ and $B^s_p$-spaces); see e.g.\ the overview in
the introduction to \cite{G90}. For $p=2$, the two scales are
identical, but for $p\ne 2$ they are related by strict inclusions:$$
H^s_p\subset B^s_p\text{ when }p>2,\quad H^s_p\supset B^s_p\text{ when }p<2.\tag 1.6
$$

When working
with operators of possibly noninteger order, it is much preferable to use
the different notations for the two scales, rather than formulating results in the
$W^{s,p}$-scale, where the definition changes when $s$ changes between integer
and noninteger values, so that mapping properties risk not being optimal.

For any open subset $\Omega $ of $\R^n$, one can define the local
variants:
$$\aligned
 H^{s}_{p,\operatorname{loc}}(\Omega )&=\{u\in\D'(\Omega )\mid \psi
 u\in  H^{s}_{p}(\R^n)\text{ for all }\psi \in C_0^\infty (\Omega
 )\},\\
 H^{s}_{p,\operatorname{comp}}(\Omega )&=\{u\in H^{s}_{p}(\R^n)\mid
 \supp u\text{ compact }\subset\Omega \},
\endaligned\tag1.7$$
and similar spaces with $B$.

In the Appendix we list anisotropic variants  of the Bessel-potential
and Besov scales with weights $(d,1)$, and explain their relation to Sobolev scales; this can
also be read as a supplementing information on the isotropic case
where $d=1$. 

Further notation from \cite{G15} is recalled in the start of Section 4
below where it is needed.

\head 2. Anisotropic symbols \endhead

When $P$ is a pseudodifferential operator on ${\Bbb R}^n$ of order $d\in{\Bbb R}$, it is a
well-known fact that if $P$ is elliptic, the solutions $u$ to the
equation $Pu(x)=g(x)$ on an open subset $\Omega $ are $d$ values more regular
than $g$, e.g., $g\in H^s_{p,\operatorname{loc}}(\Omega )$ implies
$u\in H^{s+d}_{p,\operatorname{loc}}(\Omega )$ for $s\in{\Bbb R}$. 
(Cf.\ e.g.\ Seeley \cite{S65}, Kohn and Nirenberg 
\cite{KN65},
H\"ormander \cite{H66,H85} and the exposition in Taylor \cite{T81}). It
was one of the purposes of setting up the rules of symbol calculus
to have easy access to such
regularity results. 

For the parabolic (heat operator)  problem $Pu(x,t)+\partial_tu(x,t)=f(x,t)$ on ${\Bbb R}^{n+1}$ it is
not quite as well-known what there holds of regularity. In the 
differential operator case,  when $P$ is strongly elliptic, then
$P+\partial_t$ and its solution operator belong
to a natural class  of anisotropic $\psi $do's on ${\Bbb R}^{n+1}$
where there are straightforward results.
But if $P$ is truly pseudodifferential, the symbol of $P+\partial_t$
does not satisfy all the estimates required in a standard $\psi $do
calculus on ${\Bbb R}^{n+1}$, but
something weaker (see the discussion in \cite{G86,G96} in Remark 1.5.1ff
and at the end of Section 4.1, with references). 
The operators were analyzed briefly in an $L_2$-framework in
\cite{G86,G96},
Sect.\ 4.2 (which focused on kernel estimates). 
 The mapping properties were extended to $L_p$-based spaces
 $H^{(s,s/d)}_p$ and $B^{(s,s/d)}_p$ in
 \cite{G95} (cf.\ Th.\ 3.1(1) there), for operators $P$ of integer order $d$,
in connection with a study of boundary value problems in the Boutet
 de Monvel framework. 
This depended on a symbol calculus carried over from the calculus
developed in the book \cite{G86,G96}, and relied on $L_p$-boundedness
theorems of Lizorkin \cite{L67} and Yamazaki \cite{Y84}.

In the present paper we are interested in heat problems with operators $P$ of primarily noninteger
order (such as $(-\Delta )^a$), not covered by \cite{G95}. 
There is the question of mapping properties of $P+\partial_t$, and of the
existence of (approximate) solution operators under suitable
parabolicity hypotheses, that allow showing regularity of solutions. 
In addition there is the question of {\it local regularity}, often
shown by use of commutations with cut-off functions. All this can be
handled by setting up a
systematic calculus, including composition rules.
Let us now present the appropriate symbols and
estimates.

In the following, $d\in\rp$ is fixed. The basic anisotropic invertible symbol in
the calculus is $\{\xi ,\tau \}$, with definition
$$
\{\xi ,\tau \}\equiv (\ang\xi ^{2d}+\tau ^2)^{1/(2d)},\tag2.1
$$
leading to the ``order-reducing'' operators
$$
\Theta ^su=\OP(\{\xi ,\tau \}^s)u\equiv \F^{-1}_{(\xi ,\tau )\to
(x,t)}(\{\xi ,\tau \}^s\F_{(x,t)\to (\xi ,\tau )}u),
\tag2.2$$
for all $s\in{\Bbb R}$. Then the anisotropic Bessel-potential spaces
can be defined by$$
H^{(s,s/d)}_p(\Rn\!\times\!\Bbb R)=\Theta ^{-s}L_p
(\Bbb R^{n+1});\tag2.3
$$
see more about such spaces and the related Besov family
$B^{(s,s/d)}_p$ below in the Appendix.

\proclaim{Definition 2.1} Let $m$ and $\nu \in {\Bbb R}$. The space $S^{m,\nu }_{1,0}({\Bbb R}^{2n+1})$
of $d$-anisotropic (or just anisotropic) symbols of order $m$ and with 
{\bf regularity number} $\nu $ consists of the $C^\infty $-functions $h(x,\xi ,\tau
)$ on ${\Bbb R}^{n}\times{\Bbb R}^n\times{\Bbb R}$ satisfying the estimates
$$
|\partial_x^\beta \partial_\xi ^\alpha \partial _\tau^jh(x,\xi ,\tau
 )|\le C_{\alpha ,\beta ,j}(\ang\xi ^{\nu -|\alpha |}+\{\xi ,\tau \}^{\nu -|\alpha |})\{\xi ,\tau \}^{m-\nu -dj},    \tag2.4
$$
for all indices $\alpha ,\beta \in{\Bbb N}_0^n$, $j\in{\Bbb N}_0$.

In particular, in case $m=0$ and $\nu \ge 0$, the symbols satisfy
$$
|\partial_x^\beta \partial_\xi ^\alpha \partial _\tau^jh(x,\xi ,\tau
 )|\le C_{\alpha ,\beta ,j}\ang\xi ^{ -|a|}\{\xi ,\tau \}^{ -dj}\le C_{\alpha ,\beta ,j}\ang\xi ^{ -|a|}\ang\tau ^{ -j}.    \tag2.5
$$\endproclaim

These symbol classes are very similar to those introduced in
\cite{G86,G96}, Section 2.1. The difference is that  $\ang{\xi ,\mu
}$ there has been replaced by  the anisotropic $\{\xi ,\tau \}$ here, and that the rule for differentiation
in $\tau $ is that it lowers the order by $dj$ instead of $j$
(still without changing the regularity number).
Therefore the symbols allow a very similar calculus. 

Let us first show that some symbols of interest lie in these classes.

\proclaim{Lemma 2.2} 

$1^\circ$ For $s\in{\Bbb R}$, the estimates {\rm (2.4)} are
satisfied by $\{\xi ,\tau \}^{s}  $ with $\nu =2 d$, $m=s$.

$2^\circ$ Let $p(x,\xi )$ be a $\psi $do symbol in $S^m_{1,0}({\Bbb
R}^{2n})$ for some $m\in{\Bbb R}$. Then, considered as a symbol on
${\Bbb R}^{2n+1}$, constant in $\tau $, it belongs to  $S^{m,m}_{1,0}({\Bbb
R}^{2n+1})$. 

$3^\circ$ Let $p(x,\xi )$ be a $\psi $do symbol in $S^d_{1,0}({\Bbb
R}^{2n})$. 
Then
$h(x,\xi ,\tau )=p(x,\xi )+i\tau $ satisfies the estimates $(2.4)$ with $\nu =d$,
$m=d$, i.e., belongs to $S^{d,d}_{1,0}({\Bbb
R}^{2n+1})$. 

Moreover, if $|p(x,\xi )+i\tau |\ge c(\ang\xi ^d+|\tau |)$ with $c>0$, 
then $(p(x,\xi )+i\tau
)^{-1}$ satisfies the estimates with $\nu =d$, $m=-d$.
\endproclaim

\demo{Proof} $1^\circ$. Denote for short
$$
\ang\xi =\sigma ,\quad \{\xi ,\tau \}=\kappa ,\tag2.6
$$
and observe that
$$
(\sigma ^{\nu -|\alpha |}+\kappa ^{\nu -|\alpha |})\kappa ^{m-\nu }\simeq \cases \kappa
^{m-|\alpha |}\text{ if }\nu \ge 0,\\
\sigma ^{\nu -|\alpha |}\kappa
^{m-\nu }\text{ if }\nu \le 0.\endcases
\tag2.7
$$

For $\{\xi ,\tau \}=(\ang\xi ^{2d}+\tau ^2)^{1/2d}$ itself, we have that
$$
\aligned
\partial_{\xi _j}(\ang\xi ^{2d}+\tau ^2)^{1/2d}&=\tfrac 1{2d}(\ang\xi
^{2d}+\tau ^2)^{(1/2d)-1}2d\ang\xi ^{2d-1}\partial _{\xi _j}\ang\xi
=\kappa ^{1-2d}\sigma ^{2d-1}\partial _{\xi _j}\ang\xi,\\
\partial_{\tau }(\ang\xi ^{2d}+\tau ^2)^{1/2d}&=\tfrac 1{2d}(\ang\xi
^{2d}+\tau ^2)^{(1/2d)-1}2\tau =c\kappa ^{1-2d}\tau ,
\endaligned
$$
which satisfy the estimates with $m=1$ and $\nu =2d$ (note that $\partial_{\xi _j}\ang\xi $ is
bounded and $\ang\tau \le \kappa $). Further
differentiations give linear combinations of such expressions, where
$\partial _\xi ^\alpha $ produces a factor $\sigma ^{2d-|\alpha |}$ in
at least one of the terms, whereas $\partial_\tau ^j$ only gives factors comparable
with powers of $\kappa $. The weakest term in each expression is the
one with the lowest power of $\sigma $; for $\partial_\xi
^\alpha \partial_\tau ^j\kappa $ the power is $2d-|\alpha |$, so (2.4) holds with
$m=1$, $\nu =2d$.

For $\kappa ^s$ we use that its derivatives are linear combinations of products of
$\kappa ^{s-k}$ ($k\in{\Bbb N}_0$) with expressions as above.

$2^\circ$. We have for all $\alpha $:
$$
|\partial_\xi ^\alpha  p |\le C\ang\xi ^{m-|\alpha |}\le C(\sigma ^{m-|\alpha |}+\kappa ^{m-|\alpha |})\kappa ^{m-m},\tag2.8
$$
showing the asserted estimates. 

$3^\circ$.
For $p(x,\xi )+i\tau $ we have that clearly $|p+i\tau |\le
C\kappa ^d$. Moreover, $\partial_\xi ^\alpha  (p+i\tau )$ satisfies
the estimates (2.8) with $m=d$ for $|\alpha |>0$, and 
$
\partial_\tau (p+i\tau )=i$, the higher derivatives being 0. This
shows the first assertion.

For the second assertion, the
given hypothesis shows that $|(p(x,\xi )+i\tau )^{-1}|\le C\kappa
^{-d}$;
then since  the derivatives produce negative integer powers $(p(x,\xi )+i\tau )^{-k}$ times
derivatives of $p(x,\xi )+i\tau $, the assertion is seen using the first
 estimates.
\qed
\enddemo

We have as in \cite{G86,G96}, Prop.\ 2.1.5:

\proclaim{Lemma 2.3}
The product
of two symbols of order and regularity $m,\nu $ resp.\ $m',\nu'
$ is of order $m''=m+m'$ and regularity $\nu ''=\min\{\nu ,\nu ',\nu
+\nu '\}$.

Likewise, when $h(x,\xi ,\tau )\in S^{m,\nu }_{1,0}({\Bbb
R}^{2n+1})$ and  $h'(x,\xi ,\tau )\in S^{m',\nu '}_{1,0}({\Bbb
R}^{2n+1})$, then
the Leibniz product
$$
h''(x,\xi ,\tau )\sim  {\sum}_{\alpha \in{\Bbb N}_0^n}\tfrac1{\alpha !}D_\xi ^{\alpha }h(x,\xi ,\tau )\partial_x ^{\alpha }h'(x,\xi ,\tau )\tag2.9
$$
belongs to $S^{m'',\nu ''}_{1,0}({\Bbb
R}^{2n+1})$.
\endproclaim

\demo{Proof}
The first assertion is seen by use of the elementary fact that
$$
((\sigma /\kappa )^\nu +1)((\sigma /\kappa )^{\nu'} +1)\le 3(\sigma /\kappa )^{\nu ''}+1.\tag2.10
$$
(Details in the proof of \cite{G86,G96}, Prop.\ 2.1.6.) The second
assertion now follows by application of the first assertion to each
term in the asymptotic series.\qed
\enddemo 

Here the Leibniz product represents the symbol of the composition of 
$\OP(h)$ and $\OP(h')$, as in \cite{G86,G96},
(2.1.56). Note that
since the symbols are constant in $t$, there are only terms with
$x,\xi $-derivatives. 

\example{Remark 2.4} The proofs in \cite{G96} are formulated in the framework of globally estimated $\psi $do's
of H\"ormander \cite{H85}, Sect.\ 18.1 (estimates in $x\in{\Bbb R}^n$), whereas the proofs in \cite{G86} are based on a more pedestrian local $\psi $do
calulus.
The global calculus has the advantage that remainders of order
$-\infty $ are 
treated in a simpler way, and the
Leibniz product has a precise meaning. We will take advantage of this
fact in the following, and refer to \cite{H85} and \cite{GK93,G96} for
more detailed explanations.
\endexample

\example
{Remark 2.5}
A classical $\psi $do $P$ of order $d>0$ is said to be {\it strongly
elliptic}, when the principal symbol $p_0(x,\xi )$ takes values in a
sector $V_\delta =\{z\in{\Bbb C}\mid |\operatorname{arg}z|\le \frac\pi
2-\delta \}$ for $|\xi |\ge 1$, some $0<\delta <\frac\pi 2$;
equivalently, $\operatorname{Re}p_0(x,\xi )\ge c_0|\xi |^d$ for $|\xi
|\ge 1$,
with $c_0>0$. (Recall that we take $p_0$ to be homogeneous in $\xi $
for $|\xi |\ge 1$, and $C^\infty $.) It is not hard to choose the
extension of the homogeneous function into $|\xi |\le 1$ to keep
satisfying $p_0(x,\xi )\in V_\delta $, with
$\operatorname{min}_{|\xi |\le 1}\operatorname{Re}p_0(x,\xi )>0 $. Then
$p_0(x,\xi )+i\tau $ satisfies $|p_0(x,\xi )+i\tau |\ge c(\ang\xi
^d+|\tau |)$ with $c>0$. The operators $P+\partial_t$ and symbols
$p+i\tau $ are called {\it parabolic} in this case. See also the discussion in
\cite{G86,G96}, Definition 1.5.3 ff.
\endexample

In the parabolic case there is a parametrix symbol (symbol of an
approximate inverse):

\proclaim{Lemma 2.6}
Let  $h(x,\xi ,\tau )=p(x,\xi )+i\tau $, where $p(x,\xi )$ is a
classical strongly elliptic $\psi $do symbol of order $d\in\rp$ on ${\Bbb R}^n$.Then  there is
a parametrix symbol
$k(x,\xi ,\tau )$ such that  $(p+i\tau )\#k(x,\xi ,\tau )-1$ and $ k(x,\xi
,\tau )\#(p+i\tau )-1$ are in $ \bigcap_{k\in {\Bbb
N}_0}S^{-k,d-k}_{1,0}({\Bbb R}^{2n+1})$ (i.e., they are $0$ modulo
regularity $d$). 
\endproclaim

\demo{Proof} As noted in Remark 2.5, $p$ has a principal symbol $p_0$
that is nonvanishing and takes values in a closed subsector of
$\{\operatorname{Re}z>0\}$, whereby $p_0+i\tau $ and its inverse $(p_0+i\tau )^{-1}$ are as in Lemma 2.2
$3^\circ$. Here $k_0=(p_0+i\tau )^{-1}$ is the principal term in a parametrix symbol
$k$ for $h$. The construction of the remaining terms in a full
parametrix symbol is a  standard
construction
similar to
the proof of \cite{G86,G96}, Th.\ 2.1.22.\qed     
\enddemo

Now we shall account for the mapping properties of such operators
between anisotropic Bessel-potential spaces (2.3). For isotropic standard
$\psi $do's this relies on Mihlin's multiplier theorem, but for the
present operators we need instead a  Lizorkin-type multiplier
theorem allowing separate estimates in different groups of coordinates; here we shall
use the criterion of Yamazaki \cite{Y84} (see the account of the
various criteria in \cite{GK93}, Sect.\ 1.3):

\proclaim{Lemma 2.7}(\cite{Y84}) Let $n'\in{\Bbb N}$. When $a(y,\eta )$ on ${\Bbb R}^{n'}\times
{\Bbb R}^{n'} $ satisfies 
$$
|\partial _y^\beta\eta _j^{a_j}\partial_{\eta  _j}^{a_j}a(y,\eta )|\le
 C_{\beta ,j}, \quad j=1,\dots,n', \;\alpha _j\le n'+1,\; |\beta |\le 1,\tag2.11  
$$
then $\OP(a)$ is bounded in $L_p(\R^{n'})$ for $1<p<\infty $.
\endproclaim

It is a space-dependent variant of Lizorkin's criterion.
It will be used with $n'=n+1$.

\proclaim{Theorem 2.8} Let $h(x,\xi ,\tau )\in S^{m,\nu }_{1,0}({\Bbb
R}^{2n+1})$ with $\nu \ge 0$. Then $H=\OP(h(x,\xi, \tau ))$ is continuous:
$$
\OP(h(x,\xi, \tau ))\colon  H^{(s,s/d)}_p(\Rn\!\times\!\Bbb R)\to
H^{(s-m,(s-m)/d)}_p(\Rn\!\times\!\Bbb R),\text{ for all }s\in{\Bbb R}.\tag2.12
$$

\endproclaim

\demo{Proof} By Lemma 2.3, $$
\{\xi ,\tau \}^{s-m}\#h(x,\xi ,\tau
)\#\{\xi ,\tau \}^{-s}\in S^{0,\min\{\nu ,2d\} }_{1,0}({\Bbb
R}^{2n+1}).\tag2.13
$$ 
In view of (2.5), it satisfies (2.11) with $n'=n+1$,
$y=(x,t)$ and $\eta =(\xi ,\tau )$, hence
$$
H_1=\Theta ^{s-m}H\Theta ^{-s}
$$
is bounded in $L_p({\Bbb R}^{n+1})$. By the definition of the spaces
in (2.3), it follows  that (2.12) holds.\qed
\enddemo

The proof is of course simpler for $x$-independent operators, where
(2.13) is just a product.

As a corollary, we have:

\proclaim{Corollary 2.9} Under the hypotheses in Theorem {\rm 2.8},
$H$ also maps
$$
\OP(h(x,\xi, \tau ))\colon  B^{(s,s/d)}_p(\Rn\!\times\!\Bbb R)\to
B^{(s-m,(s-m)/d)}_p(\Rn\!\times\!\Bbb R),\text{ for all }s\in{\Bbb R}.\tag2.14
$$
\endproclaim

\demo{Proof} This follows from (2.12) by use of real interpolation as
in the third line of (A.7). \qed
\enddemo

\example {Example 2.10} Consider the operator $H_0=(-\Delta
)^a+\partial_t=\OP(|\xi |^{2a}+i\tau )$. Introducing the smooth
positive modification
$[\xi ]$ of $|\xi |$:
$$
[\xi ]\text{ is }C^\infty \text{ and }\ge 1/2 \text{ on }{\Bbb R}^n,\;
[\xi ]=|\xi |\text{ for }|\xi |\ge 1,\tag 2.15
$$
and setting $r(\xi )=|\xi |^{2a} -[\xi ]^{2a}$ (supported for $|\xi
|\le 1$) we can write $H_0$ as 
$$
H_0=H_0'+\Cal R,\quad H_0'=\OP([\xi ]^{2a}+i\tau ),\; \Cal R=\OP(r(\xi )).\tag 2.16
$$
Here $H_0'$ satisfies  $3^\circ$ of Lemma 2.2 with $d=2a$,
and $\Cal R$ is smoothing in ${\Bbb R}^n$, since its symbol is $O(\ang\xi
^{-N})$ for all N. Then 
$$
H_0'\colon  H^{(s,s/(2a))}_p(\Rn\!\times\!\Bbb R)\to
H^{(s-2a,s/(2a)-1)}_p(\Rn\!\times\!\Bbb R),\text{ for all }s\in{\Bbb R},\tag2.17
$$
by Theorem 2.8. For $\Cal R$, we note that
$$
|\{\xi ,\tau \}^{s-2a}r(\xi )\hat u(\xi ,\tau )|\le C|\{\xi ,\tau \}^{s}\hat u(\xi ,\tau )|,
$$
 since $2a>0$, so $\Cal R$ also has the continuity in (2.17) for all
 $s$. It follows that
$$
H_0\colon  H^{(s,s/(2a))}_p(\Rn\!\times\!\Bbb R)\to
H^{(s-2a,s/(2a)-1)}_p(\Rn\!\times\!\Bbb R),\text{ for all }s\in{\Bbb R}.\tag2.18
$$
There is the parametrix $K_0'=\OP(([\xi ]^{2a}+i\tau )^{-1})$; it clearly maps
continuously in the opposite direction of (2.18). 

The mapping
properties also hold with $H$-spaces replaced by $B$-spaces
throughout.

A similar proof works for  symbols $p(\xi )$ that are positive for
$\xi \ne 0$ and homogeneous of degree $2a$.
\endexample

\example{Example 2.11} Here are some other simple examples, to which the theory
applies: $P=(-\Delta +b(x)\cdot \nabla+c(x))^a$ of order $2a$
(fractional powers of a perturbed Laplacian), and $P=(-\Delta
)^\frac12+b(x)\cdot \nabla$ of order 1
(the square-root Laplacian with drift). The coefficients
$b(x)=(b_1(x),\dots,b_n(x))$ and $c(x)$ are taken smooth, real and bounded
with bounded derivatives. The symbol $|\xi |+ib(x)\cdot \xi $ is
complex, with real part %$\ge $
$|\xi |$, hence strongly elliptic, and  %\linebreak
$|[\xi ]+ib(x)\cdot \xi +i\tau |\ge c(\ang\xi +|\tau |)$
with $c>0$. Another example is $P=-\Delta +(-\Delta )^{\frac12}$; here $d=2$.
\endexample

\example{Remark 2.12}
There is an important work of Yamazaki \cite{Y86} dealing with
quasi-ho\-mo\-ge\-ne\-ous $\psi $do's (and para-differential generalizations),
acting in associated quasi-homogeneous variants of the Triebel-Lizorkin spaces $F^s_{p,q}$ and Besov
spaces $B^s_{p,q}$ with  $0< p,q\le \infty $; the spaces are defined by refined
techniques involving dyadic decompositions in $\xi $-space (see also
Triebel \cite{T78} for accounts of function spaces, and
Schmeisser and Triebel \cite{ST87} for anisotropic variants). However,
it seems that the $\psi $do's in \cite{Y86} do not include our cases, since
the symbol spaces require that high  derivatives of the symbols are
dominated by $\ang{(\xi ,\tau )}^{-N}$ where any $N$ is reached, in
contrast to our estimates (2.4). (Cf.\ \cite{Y86}, Definition p.\ 157--158.)

\endexample

Much of what is said above could also be  carried through in cases where $P$
is allowed to depend moreover on $t$ (allowing a treatment of nonstationary
heat equations); in particular, Lemma 2.6 could easily be generalized. We
have kept $P$ stationary here in order to draw directly on the proofs
in \cite{G96}, and leave the nonstationary case for future investigations.

\head 3. Parabolic regularity \endhead

Now we can show regularity and local regularity of solutions to
parabolic heat equations.

\proclaim{Theorem 3.1} Consider a classical strongly elliptic $\psi $do  $P=\OP(p(x,\xi
))$ of order $d\in\rp$ on ${\Bbb R}^n$.
 Let $s\in{\Bbb R}$. 

$1^\circ$ If
$u\in H^{(s,s/d)}_p(\Rn\!\times\!\Bbb R) $ satisfies
$$
(P+\partial_t)u=f\text{ with }f\in H^{(s,s/d)}_p(\Rn\!\times\!\Bbb R),\tag3.1
$$
 then 
$$u\in H^{(s+d,s/d+1)}_p(\Rn\!\times\!\Bbb R).\tag3.2
$$

$2^\circ$ The implication from {\rm (3.1)} to {\rm (3.2)} also holds if
$u$ merely satisfies $u\in
H^{(r,r/d)}_p(\Rn\!\times\!\Bbb R) $ for some large negative $r$ (this
holds in
particular if $u\in\E'({\Bbb R}^{n+1})$ or e.g.\ $u\in L_p({\Bbb R}, \E'({\Bbb R}^n))$).

Similar statements hold with $H_p$ replaced by $B_p$ throughout.
\endproclaim 

\demo{Proof} 
Let $k(x,\xi ,\tau )$ be
a parametrix symbol according to Lemma 2.6. Then $H=P+\partial_t$ and
$K=\OP(k(x,\xi ,\tau ))$ satisfy
$$
KH=I+\Cal R_1, \text{ where }\Cal R_1=\OP(r_1),\; r_1(x,\xi ,\tau )\in
{\bigcap }_kS^{-k,d-k}_{1,0}({\Bbb R}^{2n+1}). 
$$
With the notation (2.6), we have that
$r_1(x,\xi ,\tau )$ satisfies the estimates 
$$
\multline
|\partial_\xi ^{\alpha }\partial_\tau ^jr_1|\le C (\sigma ^{d-k-|\alpha |} +\kappa  ^{d-k-|\alpha |})\kappa
 ^{-k-(d-k)-dj}
=C (\sigma ^{d-k-|\alpha |}\kappa ^{-d-dj} +\kappa  ^{-k-|\alpha
|-d-dj})\\
\le C '\sigma ^{d-k-|\alpha |}\kappa ^{-d-dj}
\le C'\sigma ^{-|\alpha |}\kappa ^{-d-dj} \text{ when }k\ge d,.\endmultline
$$
Hence $r_1\in S^{-d,0}_{1,0}({\Bbb R}^{2n+1})$ (the lowest    order we
can assign with a nonnegative regularity number). %in particular that 

Now if $u $ and $f$ are given as in $1^\circ$, then
$$
KHu=u+\Cal R_1u,
$$
where $KHu=Kf\in  H^{(s+d,s/d+1)}_p(\Rn\!\times\!\Bbb R)$ since $K$ is
of order $-d$ with regularity number $d$, and $\Cal R_1u\in
H^{(s+d,s/d+1)}_p(\Rn\!\times\!\Bbb R)$ since $\Cal R_1$ is
of order $-d$ with regularity number $0$, cf.\ Theorem 2.8. It follows that
$u\in  H^{(s+d,s/d+1)}_p(\Rn\!\times\!\Bbb R)$. This shows $1^\circ$.

$2^\circ$. Recall that a distribution in $\E'({\Bbb R}^{n+1})$, i.e.,  with compact support in ${\Bbb
R}^{n+1}$, is of
finite order, in particular lies in $H_p^{-M}(\R^{n+1})$ for some
 $M\ge 0$. A distribution in $\E'({\Bbb R}^n)$ lies in
 $H_p^{-M}({\Bbb R}^n)$ for some $M\ge 0$.  Note moreover that
$$
H_p^{-M}(\R^{n+1})\subset H_p^{(r,r/d)}(\Rn\!\times\!\Bbb R),
$$
with $r=-M$ if $d\le 1$, $r=-M/d$ if $d\ge 1$. Also, $L_p({\Bbb
 R},H_p^{-M}(\Rn))\subset H_p^{(-M,-M/d)}(\Rn\!\times\!\Bbb R)$. So we can assume
 $u\in H_p^{(r,r/d)}(\Rn\!\times\!\Bbb R)$. 

We use a bootstrap method, iterating applications
of $1^\circ$, as follows: If $r\ge s$, the statement is covered by
$1^\circ$. If $r<s$, we observe that a fortiori $f\in
H^{(r,r/d)}_p(\Rn\!\times\!\Bbb R) $. An application of $1^\circ$ then
gives the conclusion  $u\in H^{(r+d,r/d+1)}_p(\Rn\!\times\!\Bbb
R)$. Here if $r_1=r+d\ge s$, we need only apply $1^\circ$ to reach the
desired conclusion. If $r_1<s$, we repeat the argument, concluding that
 $u\in H^{(r_2,r_2/d)}_p(\Rn\!\times\!\Bbb R)$ for $r_2=r+2d$. The
 argument is repeated until $r_k=r+kd\ge s$. 

The proofs in the scale of $B_p$-spaces follow by replacing $H_p$ by
$B_p$ throughout.
\qed

\enddemo

The conclusion is best possible, in view of the forward mapping
properties in Theorem 2.8 and Corollary 2.9. 

We can also show a local regularity result.

\proclaim{Theorem 3.2} Let $P$ be as in Theorem {\rm 3.1}.
 Let $s\in{\Bbb R}$, and let $\Sigma $ be an open subset of ${\Bbb
 R}^{n+1}$.
 If
$u\in H^{(s,s/d)}_p(\Rn\!\times\!\Bbb R) $ satisfies
$$
(P+\partial_t)u|_{\Sigma }\in  H^{(s,s/d)}_{p,\operatorname{loc}}(\Sigma ),\tag 3.3
$$
 then $u|_\Sigma \in H^{(s+d,s/d+1)}_{p,\operatorname{loc}}(\Sigma )$.
\endproclaim

\demo{Proof} Notation as in Theorem 3.1 will be used. We have to show that for any $(x_0,t_0)\in\Sigma $, there
is a function  $\psi \in C_0^\infty (\Sigma )$ that is 1 on a
neighborhood of $(x_0,t_0)$ such that
$\psi u\in
H^{(s+d,s/d+1)}_{p}(\Rn\!\times\!\Bbb R )$.

Let $(x_0,t_0)\in\Sigma  $ and let $B^n_j=\{x\in{\Bbb R}^n\mid |x-x_0|<r/j\}$ and
$B^1_j=\{t\in{\Bbb R}\mid |t-t_0|< r/j\}$ for $j=1,2,\dots$, with $r>0$ so small that
$\overline{B^n_1\times B^1_1}\subset \Sigma $. For $j\in{\Bbb N}$,
define functions $\psi _j$ such that
$$
\aligned
&\psi _j(x,t)=\varphi _j(x)\varrho _j(t)\in C_0^\infty (\Sigma )\text{
with }\\
&\supp
\varphi _j\subset B^n_j,\;\varphi _j(x)=1\text{ on  }B^n_{j+1},\;
\supp\varrho  _j\subset B^1_j,\;\varrho  _j(t)=1\text{ on }B^1_{j+1}.
\endaligned\tag3.4
$$

It is given that $\psi _1u\in H^{(s,s/d)}_p(\Rn\!\times\!\Bbb R) $
and $\psi _1Hu \in H^{(s,s/d)}_p(\Rn\!\times\!\Bbb R) $. Now
$$
H\psi _1u=\psi _1Hu+[H,\psi _1]u.\tag3.5
$$
 The commutator $[H,\psi _1]=H\psi _1 - \psi _1 H$ satisfies
$$
[H,\psi _1]u=\varrho _1(t)[P,\varphi  _1(x)]u+\varphi _1(x)\partial_t\varrho  _1(t)u,
$$
where $[P,\varphi  _1]$ is a $\psi $do in $x$ of order $d-1$; in the calculus on ${\Bbb R}^{n+1}$
it counts as an operator with symbol in $S^{d-1,d-1}_{1,0}$, cf.\
Lemma 2.2 $2^\circ$.

If $d\le 1$,  $[P,\varphi  _1]$ also satisfies the estimates for an operator with symbol in
 $S^{0,0}_{1,0}$. Then both terms in the right-hand side of (3.5)
 are in $H^{(s,s/d)}_p(\Rn\!\times\!\Bbb R)$,
 and an application of $K$ to (3.5) shows that $\psi _1u\in
 H^{(s+d,s/d+1)}_p(\Rn\!\times\!\Bbb R)$, as in the proof of Theorem
 3.1.

If $d>1$, $[H,\psi _1]$ is of positive order (with symbol in
$S^{d-1,d-1}_{1,0}$), sending $u$ into\linebreak $
H^{(s-d+1,s/d-1/d+1)}_p(\Rn\!\times\!\Bbb R)$. Then we can only conclude
by application of $K$ that $\psi _1u\in  H^{(s+1,s/d+1/d)}_p(\Rn\!\times\!\Bbb R)$. Here
 we need to make an extra effort. Take $\psi _2=\varphi _2(x)\varrho
 _2(t)$ as in (3.4).
Now we can write, since $\psi _1\psi _2=\psi _2$,
$$
H\psi _2u=H\psi _1\psi _2u=\psi _2H\psi _1u+[H,\psi _2]\psi _1u=\psi _2Hu+\psi _2H(\psi _1-1)u+
[H,\psi _2]\psi _1u.\tag3.6
$$
In the final expression, the first term is in $ H^{(s,s/d)}_p(\Rn\!\times\!\Bbb R)$. In the
second term, since $\psi _2(\psi _1-1)=0$,
$$
\psi _2H(\psi _1-1)=\varphi _2(x)\varrho _2(t)P(\varphi _1(x)\varrho _1(t)-\varphi _1(x)+\varphi _1(x)-1)=\varphi _2(x)\varrho _2(t)P(\varphi _1(x)-1),
$$
where 
 $\varphi  _2P(\varphi  _1-1)$ is a $\psi $do in $x$ of order $-\infty $
 so we can regard it as an operator with symbol in
$S^{0,0}_{1,0}$. Then the term lies in $ H^{(s,s/d)}_p(\Rn\!\times\!\Bbb R)$. For
the third term, $[H,\psi _2]=\varrho _2[P,\varphi  _2]+\varphi
 _2\partial_t\varrho  _2$, where $[P,\varphi  _2]$ enters as an operator with symbol in
$S^{d-1,d-1}_{1,0}$, and when this is applied to  $\psi _1u\in
H^{(s+1,s/d+1/d)}_p(\Rn\!\times\!\Bbb R)$, we get a term in $
H^{(s+2-d,(s+2-d)/d)}_p(\Rn\!\times\!\Bbb R)$. 
If $d\le 2$, we see altogether that $H\psi _2u\in 
 H^{(s,s/d)}_p(\Rn\!\times\!\Bbb R)$, and an application of $K$ as in
 Theorem 3.1 shows
 that $\psi _2u\in 
 H^{(s+d,s/d+1)}_p(\Rn\!\times\!\Bbb R)$, ending the proof.

If $d>2$, we repeat the argument, first using that
$\psi _3=\psi _3\psi _2$, leading to the information that 
 $H\psi _3u\in 
 H^{(s,s/d)}_p(\Rn\!\times\!\Bbb R)+
 H^{(s+3-d,(s+3-d)/d)}_p(\Rn\!\times\!\Bbb R)$, next $\psi _4\psi
 _3=\psi _4$, and so on,  until $j\ge d$, so that $H\psi _ju\in
 H^{(s,s/d)}_p(\Rn\!\times\!\Bbb R)$, and we can conclude that
$\psi _ju\in 
 H^{(s+d,s/d+1)}_p(\Rn\!\times\!\Bbb R)$.\qed

\enddemo

By use of embedding theorems, we can also obtain a local regularity
result in anisotropic H\"older spaces:

\proclaim{Theorem 3.3} Let $P$ be as in Theorem {\rm
3.1.}  Let $s>0$, 
%$s\in \rp\setminus {\Bbb N}$ with $s/d\notin{\Bbb N}$, 
and let $u\in 
C^{(s,s/d)}(\Rn\!\times\!\Bbb R)\cap \E'(\Rn\!\times\!\Bbb R)$.
If, for a bounded open subset $\Omega \times
 I$,
$$
(P+\partial_t)u|_{\Omega \times I }\in  C^{(s,s/d)}_{\operatorname{loc}}(\Omega \times I),\tag 3.7
$$
 then for small $\varepsilon >0$, $u|_{\Omega \times I}\in
 C^{(s+d-\varepsilon ,(s-\varepsilon )/d+1)}_{\operatorname{loc}}(\Omega \times I
)$.
\endproclaim

\demo{Proof} 
Choose $p$ so large that $s-n/p>0$. In view of the first embedding in (A.15),  (3.7) implies that (3.3) holds for $\Sigma =\Omega \times
I$ with $s$ replaced by $s-\varepsilon _1$ for $\varepsilon _1>0$.
Moreover,  $u\in H^{(s-\varepsilon _1,(s-\varepsilon _1)/d)}_p(\Rn\!\times\!\Bbb R)$.
Then Theorem 3.2 shows that $u|_{\Omega \times I} \in H^{(s-\varepsilon
_1+d,(s-\varepsilon _1)/d+1)}_{p,\operatorname{loc}}(\Omega \times
I)$. Taking $\varepsilon _1<s-n/p$, we can
use the second embedding in (A.15) to see that $u$ is locally in
the anisotropic H\"older space  $C^{(s-n/p+d-\varepsilon
_1-\varepsilon _2,(s-n/p-\varepsilon _1-\varepsilon _2)/d+1)}_{\operatorname{loc}}(\Omega \times
I)$, for  $0<\varepsilon _2<s-n/p+d-\varepsilon _1$. Since $p$ can be
taken arbitrarily large, the statement in the theorem follows.\qed
\enddemo

For the fractional Laplacian $(-\Delta )^a$ and other related singular
integral operators,
 Fer\-nan\-dez-Real and Ros-Oton showed in  \cite{FR17} a 
result comparable to Theorem 3.3, in cases where $s$ and $s/(2a)<1$:
$$
u\in \ol C^{(s,s/(2a))}({\Bbb R}^n\times I),
f\in \ol C^{(s,s/(2a))}(\Omega \times I)\implies u\in 
C^{(s+2a,s/(2a)+1)}_{\operatorname{loc}}(\Omega \times I).\tag3.8
$$
This is sharper by avoiding the subtraction by $\varepsilon $.
 On the other hand, our result
expands the knowledge in many other directions, including that
it allows not just $s,s/(2a)<1$ but the values of $s$ up to $ \infty $, and  it allows variable-coefficient
operators, and does not require the symmetries entering in the definition of
 the fractional Laplacian, but just assumes that $P+\partial_t$ is parabolic.  
We think that a removal of $\varepsilon $ would be possible in Theorem
3.3 too --- by extending the action of our  anisotropic
 $\psi $do's (with symbols with finite regularity numbers) to
 anisotropic H\"older-Zygmund spaces.

 \cite{FR17} has in Cor.\ 3.8 and 
Rem.\ 6.4 some information on high spatial regularity of $u$ when $f$ has
high spatial regularity, assuming that
the integral operator kernel has a correspondingly high regularity.

To our knowledge, the regularity results obtained above are new in
several ways: by including all classical strongly elliptic $\psi $do's $P$
of positive real orders down to zero, and by including 
all $p\in \,]1,\infty [\,$, and all $s\in{\Bbb R}$.

Observe the special cases:

\proclaim{Corollary 3.4} Let $P$ and $u$ be as in Theorem
 {\rm 3.1} $2^\circ$.

If $u$ satisfies 
$$
(P+\partial_t)u=f\text{ with }f\in L_p(\Rn\!\times\!\Bbb R),\tag3.9
$$ 
then $u\in H_p^{(d,1)}(\Rn\!\times\!\Bbb R) 
=L_p({\Bbb R}, H_p^d({\Bbb
R}^n))\cap H_p^1({\Bbb R},L_p({\Bbb R}^n
 ))$. Conversely, $u\in
H_p^{(d,1)}(\Rn\!\times\!\Bbb R)$ implies {\rm (3.9)}. 

If $u\in L_p(\Rn\!\times\!\Bbb R) $ and $\Sigma $ is an open
subset of ${\Bbb R}^{n+1}$, then
$$
(P+\partial_t)u|_{\Sigma }\in  L_{p,\operatorname{loc}}(\Sigma )\implies u|_{\Sigma }\in
H_{p,\operatorname{loc}}^{(d,1)}(\Sigma ) .\tag3.10$$
\comment
and if $u\in {\bigcup}_{\varepsilon >0}C^\varepsilon(\Rn\!\times\!\Bbb R)
$,
$$(P+\partial_t)u|_{\Omega \times I }\in  {\bigcup}_{\varepsilon >0}C^\varepsilon
(\Omega \times I )\implies u|_{\Omega \times I }\in {\bigcap}_{0<\varepsilon
<d}C_{\operatorname{loc}}^{(d-\varepsilon ,1-\varepsilon /d)}(\Omega
\times I ) .\tag3.10
$$ 
\endcomment
\endproclaim

Note that in view of the embeddings (1.6), (A.8), we also have that (3.9) implies $u\in L_p({\Bbb R}, B_p^d({\Bbb
R}^n))\cap H_p^1({\Bbb R},L_p({\Bbb R}^n ))$ if $p\ge 2$; but this
cannot be inferred when $p<2$.
There are
similar observations for (3.10).

\head
4. A global estimate for the fractional Dirichlet heat equation
\endhead

We shall not in this paper abord the question of possible extensions
of the boundary value calculations of \cite{G95} to fractional-order
situations. We will just turn to a basic global $L_p$-result for the
fractional heat equation $Pu+\partial_tu=f$ on $\Omega \times I$,
with $P$ equal to $(-\Delta )^a$ (0.4) or the generalization (0.5) with smooth $k(y)$,
and provided with a
homogeneous Dirichlet condition. It will be obtained by  combination of a functional analysis
method put forward in \cite{BWZ17} with detailed information
from our earlier studies. 

Recall first some notation from \cite{G15}:

The following subsets of
${\Bbb R}^n$ will be considered:  
 $\rnpm=\{x\in
{\Bbb R}^n\mid x_n\gtrless 0\}$ (where $(x_1,\dots, x_{n-1})=x'$), and
 bounded $C^\infty $-subsets $\Omega $ with  boundary $\partial\Omega $, and
their complements.
Restriction from $\R^n$ to $\rnpm$ (or from
${\Bbb R}^n$ to $\Omega $ resp.\ $\complement\comega$) is denoted $r^\pm$,
 extension by zero from $\rnpm$ to $\R^n$ (or from $\Omega $ resp.\
 $\complement\comega$ to ${\Bbb R}^n$) is denoted $e^\pm$. Restriction
 from $\crnp$ or $\comega$ to $\partial\rnp$ resp.\ $\partial\Omega $
 is denoted $\gamma _0$. 

We denote by $d(x)$ a function of the form $
d(x)=\operatorname{dist}(x,\partial\Omega )$ for $x\in\Omega $, $x$ near $\partial\Omega $,
extended to a smooth positive function on $\Omega $; $d(x)=x_n$ in the
case of $\rnp$.

Along with the spaces $H^s_p({\Bbb R}^n)$ defined in (1.4), we have
the two scales of spaces associated with $\Omega $ for $s\in{\Bbb R}$:
$$
\aligned
\dot H_p^{s}(\comega)&=\{u\in H_p^{s}({\Bbb R}^n)\mid \supp u\subset
\comega \},\\
\ol H_p^{s}(\Omega)&=\{u\in \D'(\Omega )\mid u=r^+U \text{ for some }U\in
H_p^{s}(\R^n)\};
\endaligned \tag4.1
$$
here $\operatorname{supp}u$ denotes the support of $u$. The definition
is also used with $\Omega =\rnp$. In most current texts, $\ol
H_p^s(\Omega )$ is denoted $H_p^s(\Omega )$ without the overline (that
was introduced along with the notation $\dot H_p$ in \cite{H65,H85}),
but we prefer to use it, since it is
 makes the notation more
clear in formulas where
both types occur. 
We recall that $\ol H_p^s(\Omega )$ and $\dot H_{p'}^{-s}(\comega)$ are dual
spaces with respect to a sesquilinear duality extending the $L_2(\Omega )$-scalar
product; $\frac1p+\frac1{p'}=1$.

A special role in the theory is played by the {\it order-reducing
operators}. There is a simple definition of operators $\Xi _\pm^t $ on
${\Bbb R}^n$, $t\in{\Bbb R}$,
$$ 
\Xi _\pm^t =\operatorname{OP}(\chi _\pm^t),\quad \chi _\pm^t=(\ang{\xi '}\pm i\xi _n)^t ;\tag 4.2
$$
 they preserve support
in $\crnpm$, respectively. The functions
$(\ang{\xi '}\pm i\xi _n)^t $ do not satisfy all the estimates for $S^{t }_{1,0}({\Bbb
R}^n\times{\Bbb R}^n)$, but lie in a space as in Definition 2.1 with
$d=1$, $\nu =1$, $m=t$, with $(\xi ,\tau )$ replaced by $(\xi ',\xi
_n)$. There is a more refined choice $\Lambda _\pm^t $
\cite{G90, G15}, with
symbols $\lambda _\pm^t (\xi )$ that do
satisfy all the estimates for $S^{ t }_{1,0}({\Bbb
R}^n\times{\Bbb R}^n)$; here $\overline{\lambda _+^t }=\lambda _-^{t }$.
The symbols have holomorphic extensions in $\xi _n$ to the complex
halfspaces ${\Bbb C}_{\mp}=\{z\in{\Bbb C}\mid
\operatorname{Im}z\lessgtr 0\}$; it is for this reason that the operators preserve
support in $\crnpm$, respectively. Operators with that property are
called "plus" resp.\ "minus" operators. There is also a pseudodifferential definition $\Lambda
_\pm^{(t )}$ adapted to the situation of a smooth domain $\Omega
$, cf.\ \cite{G15}.

It is elementary to see by the definition of the spaces $H_p^s(\R^n)$
in terms of Fourier transformation, that the operators define homeomorphisms 
for all $s$: $
\Xi^t _\pm\colon H_p^s(\R^n) \simto H_p^{s- t
}(\R^n)$, $  
\Lambda ^t _\pm\colon H_p^s (\R^n) \simto H_p^{s- t
} (\R^n)$.
The special
interest is that the "plus"/"minus" operators also 
 define
homeomorphisms related to $\crnp$ and $\comega$, for all $s\in{\Bbb R}$: 
$
\Xi ^{t }_+\colon \dot H_p^s(\crnp )\simto
\dot H_p^{s- t }(\crnp)$, $
r^+\Xi ^{t }_{-}e^+\colon \ol H_p^s(\rnp )\simto
\ol H_p^{s- t } (\rnp )$, with similar statements  for $\Lambda ^t_\pm$,
and for $
\Lambda^{(t )}_\pm$ relative to $\Omega $.
Moreover, the operators $\Xi ^t _{+}$ and $r^+\Xi ^{t }_{-}e^+$ identify with each other's adjoints
over $\crnp$, because of the support preserving properties.
There is a
similar statement for $\Lambda ^t_+$ and  $r^+\Lambda ^t_-e^+$, and for $\Lambda ^{(t )}_+$ and $r^+\Lambda ^{(
t )}_{-}e^+$ relative to the set $\Omega $.

The special {\it $\mu $-transmission spaces} were 
introduced by
H\"ormander \cite{H65} for $p=2$, cf.\ \cite{G15}; we shall just use
them here for $\mu =a$: 
$$
\aligned
H_p^{a (s)}(\crnp)&=\Xi _+^{-a }e^+\ol H_p^{s- a
}(\rnp)=\Lambda  _+^{-a }e^+\ol H_p^{s- a
}(\rnp)
,\quad  s> a -1/p',\\
H_p^{a (s)}(\comega)&=\Lambda  _+^{(-a )}e^+\ol H_p^{s- a
}(\Omega ),\quad  s> a -1/p';
\endaligned\tag 4.3
$$
they are the appropriate solution spaces for homogeneous Dirichlet
problems for elliptic operators $P$ having the $a
$-transmission property (cf.\ \cite{G15}). Note that in (4.3), $\Xi
_+^{-a}$ is applied to functions with a jump at $x_n=0$ (when
$s>a+1/p$), this results in a singularity at $x_n=0$.

The $\psi $do $P$ can be  applied to functions in the spaces in (4.1)
when they are extended by zero to all of ${\Bbb R}^n$. This is already
understood for the spaces $\dot H^s_p(\comega)$, but should be
mentioned explicitly (by an indication with $e^+$) for the spaces $\ol H^s_p(\Omega )$. Also, when
$u\in \dot H^s_p(\comega)$ and $Pu$ is considered on $\Omega $, it is
most correct to indicate this by writing $r^+Pu$. The indications
$e^+$ and $r^+$ can be left out
as an ``abuse of notation'', when they are understood from the
context; note however the importance of $e^+$ in (4.3).

Recall from \cite{G15}, Theorems 4.4 and 5.4:

\proclaim{Theorem 4.1} Let $\Omega $ be an open bounded smooth subset
of ${\Bbb R}^n$, let  $P$ be a $\psi $do on ${\Bbb R}^n$ of the form {\rm (0.5)}
 with $k\in C^\infty (S^{n-1})$
(in particular, $P$ can be equal to $(-\Delta )^a$) for
an $a>0$.
Let $1<p<\infty $, and let $P_{\operatorname{Dir},p}$ stand
for the $L_p$-Dirichlet realization on $\Omega $, acting like $r^+P$ and with domain 
$$
D(P_{\operatorname{Dir},p})=\{u\in \dot H_p^a(\comega)\mid r^+Pu\in
 L_p(\Omega )\};
 \tag4.4
$$
the operators are consistent for different $p$. Then 
$$
D(P_{\operatorname{Dir},p})= H^{a(2a)}_p(\comega)=\Lambda
_+^{(-a)}e^+\ol H^{a}_p(\Omega ).\tag4.5
$$
It satisfies (for any $\varepsilon >0$):
$$
 H^{a(2a)}_p(\comega)\cases =\dot H_p^{2a}(\comega),\text{ if
}a<1/p,\\
\subset\dot H_p^{2a-\varepsilon }(\comega),\text{ if }a=1/p,\\
\subset e^+d^a\ol H_p^a(\Omega )+\dot H_p^{2a
}(\comega),\text{ if }a>1/p,\; a-1/p\notin{\Bbb N},\\
\subset e^+d^a\ol H_p^a(\Omega )+\dot H_p^{2a-\varepsilon
}(\comega),\text{ if }a>1/p,\; a-1/p\in{\Bbb N}.
\endcases \tag4.6
$$
More precisely, in the case $a\in \,]0,1[\,$, the functions have in local coordinates where $\Omega
$ is replaced by $\rnp$ the following structure when $a>1/p$:
$$\aligned
u=w+x_n^aK_0\varphi , %\text{ if }a>1/p,
\endaligned
\tag4.7$$
where $w$ and $\varphi $
run through $\dot H_p^{2a}(\crnp)$ and $
B_p^{a-1/p}({\Bbb R}^{n-1} )$, 
and $K_0 $ is the Poisson operator $K_0\varphi =\Cal F^{-1}_{\xi
'\to x'}[\hat \varphi (\xi ')e^+r^+e^{-\ang{\xi '}x_n}]$.
\endproclaim

\demo{Proof} When $P$ is defined by (0.5) with  a smooth $k$, it is the
$\psi $do with symbol $p(\xi )$ equal to the
Fourier transform of the  kernel function $k(y/|y|)|y
|^{-n-2a}$. (For $P=(-\Delta )^a$, $k$ equals the constant $c_{n,a}$.) 
The symbol is a function
homogeneous of degree $2a$, smooth positive for $\xi \ne 0$, and {\it
even}: $p(-\xi )=p(\xi )$. The $L_2$-Dirichlet realization  $P_{\operatorname{Dir},2}$ can be defined
variationally from the sesquilinear form
$$
Q(u,v)=\tfrac12 \int_{{\Bbb
R}^{2n}}\frac{(u(x)-u(y))(\bar v(x)-\bar v(y))k((x-y)/|x-y|)}{|x-y|^{n+2a}}
\,dxdy\tag 4.8
$$
considered for $u,v\in \dot H^a(\Omega )$. As accounted for in Ros-Oton \cite{R16}, it satisfies a Poincar\'e{}
inequality over $\Omega $
so that the selfadjoint operator in $L_2(\Omega )$
induced by the Lax-Milgram lemma is a bijection from its domain to
$L_2(\Omega )$. The explanation in \cite{R16} is formulated for real functions, but the
operator defined in $L_2(\Omega ,{\Bbb C})$ is {\it real} in the sense that it
maps real functions to real functions, and it can be retrieved from
the definition on $L_2(\Omega ,{\Bbb R})$ by linear extension.

Since $p(\xi )$ is even and homogeneous of degree $2a$, it satisfies
the condition in \cite{G15} for having the $a$-transmission property,
and since it is positive, it has factorization index $a$ (since
$a_+=a_-$ in \cite{G15}, (3.3)--(3.4)).
Then, considering its action in $L_p$-spaces, 
the description of the domain in (4.5) follows from \cite{G15}, Th.\ 4.4 with $m=2a$, $\mu
_0=a$, $s=2a$. The consistency for various $p$ holds as a general
property of pseudodifferential operators. The statement (4.6) is from \cite{G15}, Th.\ 5.4, with
$\mu =a$, $s=2a$.

 The information (4.7) is a consequence of the proof
given there; let us give a direct explanation here: It is well-known
that when
$1/p<a<1+1/p$,  the functions $v\in \ol
H_p^a(\rnp)$ have a 
first trace $\gamma _0v\in B_p^{a-1/p}({\Bbb R}^{n-1})$,
and that $v-K\gamma _0v\in \dot H_p^a(\crnp)$, when $K$ is a continuous
right inverse of $\gamma _0$. In fact, the functions $v\in \ol H_p^a(\rnp)$
are exactly the functions of the form $v=g+K\varphi $, where  $g$
runs through $\dot H_p^a(\crnp)$ and $\varphi $ runs through $
B_p^{a-1/p}({\Bbb R}^{n-1})$. Take as $K$ the Poisson operator $K_0$, $$
K_0\varphi
=\Cal F^{-1}_{\xi 
'\to x'}[\hat \varphi (\xi ')e^+r^+e^{-\ang{\xi '}x_n}]
=\Cal F^{-1}_{\xi 
\to x}[\hat \varphi (\xi ')(\ang{\xi '}+i\xi _n)^{-1}].\tag4.9$$ 
Now to
describe $ H^{a(2a)}_p(\crnp)=\Xi 
_+^{-a}e^+\ol H^{a}_p(\rnp )$, we use that 
$$
\Xi _+^{-a}K_0\varphi
=\Cal F^{-1}_{\xi 
\to x}[(\ang{\xi '}+i\xi _n)^{-a}\hat \varphi (\xi ')(\ang{\xi '}+i\xi
_n)^{-1}]
=c_ax_n^aK_0\varphi ,$$ 
cf.\ \cite{G15}, (2.5) and (5.16). Moreover, $\Xi _+^{-a}\dot
H_p^a(\crnp)=\dot H_p^{2a}(\crnp)$. This shows the representation in
(4.7). (We are using the formulas from \cite{G15} with $[\xi ']$
replaced by $\ang{\xi '}$, which works equally well, as shown there.)
\qed
\enddemo

We underline that for $a>1/p$, $D(P_{\operatorname{Dir},p})$  contains
not only $\dot H^{2a}_p(\comega)$, but also functions of the form
$e^+d^az$ with $z\in \ol H^a_p(\Omega )$, not in $ \ol H^s_p(\Omega )$
for $s>a$. (On the interior, the functions are in
$H^{2a}_{p,\operatorname{loc}}(\Omega )$, by elliptic regularity.)

The operators $P_{\operatorname{Dir},p}$ are bijective for all $p$;
for $p=2$ this holds 
since the sesquilinear form defining
$P_{\operatorname{Dir},2}$  has positive lower bound by the Poincar\'e{} inequality, 
and for the other $p$ it follows in view of \cite{G14}, Th.\ 3.5, on
the stability of kernels and cokernels when spaces change. 

Next, we will describe $P_{\operatorname{Dir},p}$ from a functional
analysis point of view. We already have the definition of 
$P_{\operatorname{Dir},2}$ from the sesquilinear form (4.8).
 The associated
quadratic form $Q(u,u)$ is denoted $Q(u)$, 
$$
Q(u)=\tfrac12\int_{{\Bbb
R}^{2n}}\frac{|u(x)-u(y)|^2k((x-y)/|x-y|)}{|x-y|^{n+2a}}\, dxdy\text{
on }\dot H^{a}_p(\comega).\tag4.10
$$
 
In view of the positivity, $-P_{\operatorname{Dir},2}$ generates a
strongly continuous semigroup $e^{-tP_{\operatorname{Dir},2}}$ of
contractions in $L_2(\Omega ,{\Bbb R})$. 

We shall now see that the form $Q(u,v)$ is a so-called Dirichlet form, as defined in Davies
\cite{D89} and Fukushima, Oshima and Takeda \cite{FOT94}. This is
observed e.g.\ in Bogdan, Burdzy and Chen \cite{BBC03} for a related
form defining the regional fractional Laplacian. It means
that $Q$ has the Markovian property (cf.\ \cite{FOT94}, pages 4--5):

\proclaim{Definition 4.2} A closed nonnegative symmetric form $E(u,v)$
with domain $D(E)\subset L_2(\Omega ,{\Bbb R})$ is said to be
Markovian, if for any $\varepsilon >0$ there exists a function
$\varphi _\varepsilon $ on ${\Bbb R}$ taking values in $[-\varepsilon
,1+\varepsilon ]$ with $\varphi _\varepsilon (t)=t$ on $[0,1]$ and
$0\le \varphi _\varepsilon (t)-\varphi _\varepsilon (s)\le t-s$ when 
$t>s$, such that
$$
u\in D(E)\implies \varphi _\varepsilon \circ u\in D(E)\text{ and
}E(\varphi _\varepsilon u,\varphi _\varepsilon u)\le E(u,u).\tag4.11
$$
 
\endproclaim

There is an equivalent definition with $\varphi _\varepsilon \circ u$
in (4.11)  replaced by  $\min\{\max\{u, 0\},1\}$.

We can choose $\varphi _\varepsilon $ to be $C^\infty $ on ${\Bbb R}$
(but not in $C_0^\infty ({\Bbb R})$ as written in \cite{BBC03}), then
it is clear that (4.11) holds for $E=Q$.

The interest of the Markovian property here is that
then the semigroup $e^{-tP_{\operatorname{Dir},2}}$ 
extends for each $p\in \,]1,\infty [\,$  to a
semigroup $T_p(t)$  that is contractive in $L_p$-norm and bounded
holomorphic (\cite{FOT94} Th.\ 1.4.1 and \cite{D89} Th.\ 1.4.1);
with infinitesimal
generators that are consistent for varying $p$. The  generator of $T_p(t)$
is in fact $-P_{\operatorname{Dir},p}$, 
since the latter is bijective from its domain to $L_p(\Omega )$, and
the operators  $-P_{\operatorname{Dir},p}$ are consistent for varying $p$.

As pointed out in \cite{BWZ17}, one can use the existence of these
extensions to $L_p$ to apply a theorem of Lamberton \cite{L87} giving information on the
heat equation solvability.

\proclaim{Theorem 4.3} Let $P$ and $\Omega $ be as in Theorem {\rm
4.1}, with $0<a<1$, and let $I=\,]0,T[\,$ for some $T>0$.  The
Dirichlet evolution problem 
$$
\aligned
Pu+\partial_tu&=f\text{ on }\Omega \times I,\\
u&=0\text{ on }(\R^n\setminus\Omega) \times  I,\\
u&=0\text{ for }t=0,
\endaligned \tag 4.12
$$
has for any $f\in L_p(\Omega\times  I)$ a unique solution $u(x,t)\in
C^0(\overline I,L_p(\Omega ))$, which satisfies:
$$
u\in L_p(I,H_p^{a(2a)}(\comega))\cap H^1_p( I,L_p(\Omega )).\tag4.13
$$
Here $H_p^{a(2a)}(\comega)$ is the domain of
$P_{\operatorname{Dir},p}$, as  described in detail in Theorem {\rm 4.1}.
\endproclaim

\demo{Proof} In (4.12) it is tacitly understood that $u$ identifies with a function
on $\R^n\times I$ vanishing for $x\in \R^n\setminus \Omega $,
in order for the $\psi $do to be defined on $u$. 

As accounted for above, the operator $P_{\operatorname{Dir},2}$ satisfies the
hypotheses for the operator $-A$ studied in \cite{L87},
namely that $A$ generates a bounded holomorphic semigroup for $p=2$,
and induces bounded holomorphic semigroups $T_p(t)$ in $L_p(\Omega )$ that are
contractions for all $p\in \,]1,\infty [\,$, and are consistent with the
case $p=2$. Then, according to
\cite{L87}, Th.\ 1, the problem
$$
(\partial_t-A)u=f \text{ on }I,\quad u|_{t=0}=0,\tag4.14
$$
has for any $f\in L_p(\Omega \times  I)$ a solution $u(x,t)\in
C^0(\overline I,L_p(\Omega ))$ such that
$$
\|\partial_tu\|_{ L_p(\Omega \times I)}+\|Au\|_{ L_p(\Omega
\times  I)}\le C\|f\|_{ L_p(\Omega \times  I)}.\tag 4.15
$$
 The bound on the first term shows that $u\in H^1_p(
 I,L_p(\Omega ))$.
Since $\|Pv\|_{L_p(\Omega )}\ge C'\|v\|_{H_p^{a(2a)}(\comega)}$ for all $v(x)\in
 D(P_{\operatorname{Dir},p})$, the bound on the second term shows that 
 $u\in L_p( I,H_p^{a(2a)}(\comega))$.

The uniqueness of the solution is accounted for e.g.\ in \cite{LPPS15}.\qed
\enddemo

The regularity in (4.13) is optimal for $f\in L_p(\Omega \times I)$.

In view of the general rules (1.6), the
theorem implies that the solution also satisfies $u\in  L_p(
I,B_p^{a(2a)}(\comega))$ when $p\ge 2$, but hits a larger space than 
$L_p(
I,B_p^{a(2a)}(\comega))$ when $p<2$. 

For $P=(-\Delta )^a$, Biccari, Warma and Zuazua  
\cite{BWZ17} used \cite{L87} to 
show semi-local variants of (4.13): $u \in  L_p(
I,W_{\operatorname{loc}}^{2a,p}(\Omega))$ with $\partial_tu\in
L_p(I\times \Omega)$ when $p\ge 2$ or $a=\frac12$; and it holds with
$W_{\operatorname{loc}}^{2a,p}(\Omega)$ replaced by $B_{p,2,\operatorname{loc}}^{2a}(\Omega)$
when $p<2$, $a\ne \frac12$ (there is an embedding $H^s_p\subset B^s_{p,2}$
for such $p$).

When $f$ has a higher regularity, we can use Theorem 3.2 to get a
local result:

\proclaim{Theorem 4.4} Let $u$ be the solution of {\rm (4.12)} defined
in Theorem {\rm 4.3}, and let
$r=2a$ if $a<1/p$, $r= a+1/p-\varepsilon $ if $a\ge
1/p$ (for some small $\varepsilon >0$).
Then $u\in \ol H_p^{(r,r/(2a))}({\Bbb R}^n \times I)$, vanishing for $x\notin
\comega $.

If $f\in  H_{p,\operatorname{loc}}^{(s,s/(2a))}(\Omega \times I)$ for
some $0<s\le r$, then 
 $u\in  H_{p,\operatorname{loc}}^{(s+2a,s/(2a)+1)}(\Omega \times I)$.
\endproclaim

\demo{Proof} Since $e^+\ol H_p^a(\Omega )=\dot H^a_p(\comega)$ for
$a<1/p$, and is contained in $\dot H_p^{1/p-\varepsilon }(\comega)$ when
$a\ge 1/p$, $
H_p^{a(2a)}(\comega)=\Lambda _+^{-a}e^+\ol H_p^a(\Omega )\subset \dot
H_p^{r}(\comega)$, where $r$ is as defined in the theorem.
Note that $r\le 2a$, and that $H_p^{a(2a)}(\comega)\subset \dot
H_p^{a+1/p-\varepsilon }(\comega)$ in any case. 

Thus, when $u$ is a solution as in Theorem 4.3, 
$u\in \ol H_p^{(r,r/(2a))}({\Bbb R}^n \times I)$, vanishing for $x\notin
\comega $. 

In view of the initial condition $u|_{t=0}=0$, the
extension by zero for $t<0$ lies in $\ol H_p^{(r,r/(2a))}({\Bbb R}^n
\times \,]-\infty ,T[\,)$.  We can moreover extend our function for $t\ge T$ to a
function in $H_{p}^{(r,r/(2a))}(\Rn\!\times\!\Bbb R)$, and denote the fully
extended function $\tilde u$. Observe
that the values of $\tilde f=(P+\partial_t)\tilde u$ are
consistent with the values of $f$ on $\Omega \times I$, since $P$ acts
only in the $x$-direction and $\partial_t$ is local.

Now if $f\in  H_{p,\operatorname{loc}}^{(s,s/(2a))}(\Omega \times I)$ 
for some $0<s\le r$,
we can apply Theorem 3.2, and conclude that $u\in  H_{p,\operatorname{loc}}^{(s+2a,s/(2a)+1)}(\Omega \times I)$.\qed 
\enddemo

For higher values of $s$, Theorem 3.2 gives that if  $u\in
H_{p}^{(s,s/(2a))}({\Bbb R}^n\times I)$ and $f\in
H_{p,\operatorname{loc}}^{(s,s/(2a))}(\Omega \times I)$, then  $u\in
H_{p,\operatorname{loc}}^{(s+2a,s/(2a)+1)}(\Omega \times I)$,  but the
global prerequisite on $u$ may not be easy to obtain. 

\example{Remark 4.5}
In comparison with the result of \cite{FR17}, Cor.\ 1.6, in anisotropic H\"older
spaces, our study in $H^{(s,s/(2a))}_p$-spaces has the advantage
that the regularity in
$t$ of the solution of the Dirichlet problem is lifted by a full step
1 in Theorem 4.3, whereas \cite{FR17} obtains a $C^{1-\varepsilon
}$-estimate in $t$.

Also for the interior regularity we observe a better lifting in
$t$-derivatives, namely that when $a<1/p$, $t$-derivatives of order 2 are
controlled in Theorem 4.4 (taking $s=2a$). Since $0<a<1$, this holds for $p$ sufficiently
close to 1. 
For the H\"older estimates in \cite{FR17}, Th.\
1.1, second $t$-derivatives are not reached (cf.\ (3.8) above).

\endexample

\example{Remark 4.6} Theorem 4.1 is also valid for $x$-dependent
strongly elliptic $\psi $do's $P$ of order $2a$ with even symbol. Then
$P_{\operatorname{Dir},2}$ is defined from the sesquilinear form
$(Pu,v)_{L_2(\Omega )}$ on $C_0^\infty (\Omega )$, extended by closure
to $\dot H^a_2(\comega)$, and the $P_{\operatorname{Dir},p}$ are
consistent with this by (4.4). A variant of Theorem 4.3 can be shown
for $p=2$ by techniques from Lions and Magenes \cite{LM68} vol.\ 2 (we shall
explain details elsewhere), but for $p\ne 2$, other methods are needed.

\endexample
\head
Appendix. Anisotropic Bessel-potential and Besov spaces.
\endhead

In this Appendix we present the appropriate anisotropic generalizations of
Bes\-sel-potential and Besov spaces. We just give a summary of the
essentially well-known facts that are needed for the parabolic
 operator $P+\partial_t$ on ${\Bbb R}^n\times {\Bbb R}$ with a
$\psi $do $P$ on ${\Bbb R}^n$ of positive order. Let
$d\in\rp$ and $p\in\,]1,\infty [\,$. This material is taken from the appendix of
\cite{G95}, with small modifications and less focus on cases where $d$
is integer.

For $m,d\in\Bbb N$, the anisotropic Sobolev spaces
$W^{(dm,m)}_p(\Rn\!\times\!\Bbb R)$ are defined by$$\aligned
&W^{(dm,m)}_p(\Rn\!\times\!\Bbb R)=L_p\bigl(\Bbb R;H^{dm}_p(\Rn)\bigr)
\cap H^{m}_p\bigl(\Bbb R;L_p(\Rn)\bigr)\\
&=\{\,u(x,t)\in L_{p}(\Bbb
R^{n+1})\mid
\F^{-1}_{x,t}(\ang{\xi }^{dm}+\ang{\tau }^m)\hat u)\in L_p(\Bbb R^{n+1})\,\}\\
&=\{\,u(x,t)\in L_{p}(\Bbb
R^{n+1})\mid
D_x^\alpha D_t^ju\in L_p(\Bbb R^{n+1})\text{ for }|\alpha |+dj\le dm\,\}.
\endaligned\tag A.1$$
For the generalization to noninteger and negative values of the
Sobolev exponents, 
we
observe that if we define $
\{\xi ,\tau \}$ and
$\Theta  ^s=\operatorname{OP}_{x,t}(\{\xi ,\tau \} ^s)$ as in (2.1), (2.2),
then
$$
W^{(dm,m)}_p(\Rn\!\times\!\Bbb R)=\Theta  ^{-dm}L_p(\Bbb R^{n+1}).\tag A.2
$$
This is seen by use of Lizorkin's criterion \cite{L67}, cf.\ (2.11) above,
applied to the operators
with symbol $(\ang{\xi }^{dm}+\ang{ \tau }^m)
\{\xi ,\tau \}  ^{-dm}$ and $\{\xi ,\tau \} ^{dm}
(\ang{\xi }^{dm}+\ang{ \tau}^m)^{-1}$.

Clearly, $\Theta ^s\Theta ^t=\Theta ^{s+t}$
for $s,t\in \Bbb R$.

More generally, one can now define with $d\in\rp$, for any $s\in\Bbb R$,  $$
H^{(s,s/d)}_p(\Rn\!\times\!\Bbb R)=\Theta ^{-s}L_p
(\Bbb R^{n+1});\tag A.3
$$
it is an anisotropic
generalization of the Bessel-potential spaces $H^s_p(\Rn)$, and
clearly \linebreak$H^{(dm,m)}_p(\Rn\!\times\!\Bbb R)=
W^{(dm,m)}_p(\Rn\!\times\!\Bbb R)$ for $m,d\in\Bbb N$. Here we
follow the notation of Schmeisser and Triebel, cf\. e.g\. \cite{ST87} and
its references,
deviating from another extensively used notation
$L^{(s,s/d)}_p(\Rn\!\times\!\Bbb R)$ 
(as in Lizorkin \cite{L66,L67},  Nikolski\u\i{} \cite{N75},
Besov, Il'in and Nikolski\u\i{} \cite{BIN78},\dots),
where the spaces are often called Liouville spaces. See also Sadosky and
Cotlar \cite{SC66}. These spaces fit
together in {\it complex interpolation} (by an anisotropic
variant of the proof for $H^s_p(\Bbb R^{n+1})$ spaces, cf.\ Calder\'on
\cite{C63}, Schechter \cite{S67}, and  Schmeisser and Triebel
\cite{ST74} Rem.\ 4):$$\multline
[H^{(s_0,s_0/d) }_p(\Rn\!\times\!\Bbb R),H^{(s_1,s_1/d) }_p(\Rn\!\times\!\Bbb R)]_\theta= 
H^{(s_2,s_2/d) }_p(\Rn\!\times\!\Bbb R),\\
\quad s_2=(1-\theta )s_0+\theta s_1,\quad\text{for }s_0,s_1\in\Bbb
R,\, \theta \in\,]0,1[\,.\endmultline\tag A.4
$$ 

Another generalization of the $W^{(dm,m)}_p$ spaces is 
the scale of anisotropic Besov spaces
$B^{(s,s/d)}_p(\Rn\!\times\!\Bbb R)$, that can be defined as
follows:
$$\aligned
\|u\|^p_{B^{(s,s/d)}_p(\Rn\times\Bbb R)}&= 
\|u\|^p_{L_p}
+ \int_{\Bbb R^{2n+2}}\Bigl(\frac{|
u(x,t)-u(y,t)|^p}{|x-y|^{n+ps}}\\
&\quad+
\frac{| u(x,t)-u(x,t')|^p}{|t-t'|^{1+ps/d}}\Bigr)\,dxdydtdt'\text{
for }s\in\,]0,1[\,;\\
B^{(s+r,(s+r)/d)}_p(\Rn\!&\times\!\Bbb R)= \Theta ^{-r}
B^{(s,s/d)}_p(\Rn\!\times\!\Bbb R) ,\text{ for }r\in\Bbb
R,\;s\in\,]0,1[\,.
\endgathered\tag A.5
$$

The norm
can also be described in terms of dyadic decompositions 
 (see e.g\. \cite{ST87}),
or, for $s>0$, by a generalization of the integral formula in (A.5)
involving higher order differences as in (1.5) (see e.g\. Besov \cite{B64}, Solonnikov
\cite{So65}, \cite{LSU68}, \cite{BIN78}).
These spaces arise from the $W^{(dm,m)}_p(\Bbb R^n\!\times\!\Bbb
R)$ spaces by {\it real interpolation} 
(cf\. Grisvard \cite{G66, I.9}), when $m,d\in\Bbb N$:
$$
\bigl(W^{(dm,m)}_p(\Rn\!\times\!\Bbb R),L_p(\Bbb R^{n+1})\bigr)_{\theta
,p}= 
B^{((1-\theta )dm,(1-\theta )m) }_p(\Rn\!\times\!\Bbb
R);\;\theta \in\,]0,1[\,.\tag A.6
$$
Moreover, one has 
for all $d\in\rp$, all $s_0,s_1,s\in\Bbb R$ with $s_0\neq s_1$, 
all $p_0,p_1,p\in\,]1,\infty [\,$ with  $p_0\ne p_1$, all
$ \theta \in\,]0,1[\,$, setting
$s_2=(1-\theta )s_0+\theta s_1$, 
$\frac1{p_2}=\frac{1-\theta }{p_0}+\frac{\theta }{p_1}$ 
(cf\. \cite{G66} and \cite{ST87, 3.2}):
$$\aligned
\bigl(B^{(s_0,s_0/d) }_p(\Rn\!\times\!\Bbb R),B^{(s_1,s_1/d) }_p(\Rn\!\times\!\Bbb R)\bigr)_{\theta
,p}&= 
B^{(s_2,s_2/d) }_p(\Rn\!\times\!\Bbb R),\\
[B^{(s_0,s_0/d) }_p(\Rn\!\times\!\Bbb R),B^{(s_1,s_1/d) }_p
(\Rn\!\times\!\Bbb R)]_\theta&= 
B^{(s_2,s_2/d) }_p(\Rn\!\times\!\Bbb R),\\
\bigl(H^{(s_0,s_0/d) }_p(\Rn\!\times\!\Bbb R),H^{(s_1,s_1/d) }_p(\Rn\!\times\!\Bbb R)\bigr)_{\theta
,p}&= 
B^{(s_2,s_2/d) }_p(\Rn\!\times\!\Bbb R),\\
\bigl(H^{(s,s/d) }_{p_0}(\Rn\!\times\!\Bbb R),H^{(s,s/d)
}_{p_1}(\Rn\!\times\!\Bbb R)\bigr)_{\theta 
,p_2}&= 
H^{(s,s/d) }_{p_2}(\Rn\!\times\!\Bbb R).
\endaligned\tag A.7
$$

The Bessel-potential spaces and Besov spaces are interrelated by
$$
\aligned
B^{(s,s/d)}_{p}(\Rn\!\times\!\Bbb R)&\subset
H^{(s,s/d)}_{p}(\Rn\!\times\!\Bbb R) \subset B^{(s-\varepsilon
,(s-\varepsilon )/d)}_{p}(\Rn\!\times\!\Bbb R)\text{ for }p\le 2;\\
H^{(s,s/d)}_{p}(\Rn\!\times\!\Bbb R)&\subset
B^{(s,s/d)}_{p}(\Rn\!\times\!\Bbb R) \subset H^{(s-\varepsilon
,(s-\varepsilon )/d)}_{p}(\Rn\!\times\!\Bbb R)\text{ for }p\ge 2;
\endaligned\tag A.8
$$
for $s\in\Bbb R$, any $\varepsilon >0$; and $B^{(s,s/d)}_p\ne
H^{(s,s/d)}_p$ when $p\ne 2$.
 
For both types, one has the identification of {\it dual spaces} (with
$\frac 1p+\frac1{p'}=1$ as usual):
$$
\aligned
H^{(s,s/d)}_p(\Rn\!\times\!\Bbb
R)^*&\simeq H^{(-s,-s/d)}_{p'}(\Rn\!\times\!\Bbb R) \quad\text{and}\\
B^{(s,s/d)}_p(\Rn\!\times\!\Bbb
R)^*&\simeq B^{(-s,-s/d)}_{p'}(\Rn\!\times\!\Bbb R), \quad\text{for
}s\in\Bbb R.\endaligned\tag A.9
$$

For positive $s$, the spaces can moreover be described in the following way:
$$\alignedat2
\text{\rm (i)}&&\quad
H^{(s,s/d)}_{p}(\Rn\!\times\!\Bbb R)&=L_p\bigl(\Bbb R; H^{s}_p(\Bbb
R^n)\bigr)\cap H^{s/d}_p\bigl(\Bbb R;L_p(\Bbb R^n)\bigr),\; s\ge 0;\\
\text{\rm (ii)}&&\quad B^{(s,s/d)}_{p}(\Rn\!\times\!\Bbb
R)&=L_p\bigl(\Bbb R; B^{s}_p(\Bbb 
R^n)\bigr)\cap B^{s/d}_p\bigl(\Bbb R;L_p(\Bbb R^n)\bigr),\; s>0;
\endalignedat\tag A.10 
$$
cf\. Grisvard \cite{G66,G69}. The anisotropic Sobolev-Slobodetski\u\i{} spaces
$W^{(s,s/d)}_p(\Rn\!\times\!\Bbb R)$ are defined from the i\-so\-tro\-pic
ones by (cf\. e.g\. Solonnikov \cite{S87})
$$\aligned
W^{(s,s/d)}_{p}(\Rn\!\times\!\Bbb R)&=L_p\bigl(\Bbb R; W^{s}_p(\Bbb
R^n)\bigr)\cap W^{s/d}_p\bigl(\Bbb R;L_p(\Bbb R^n)\bigr),\; s\ge 0;\\
\text{thus }W^{(s,s/d)}_{p}&=H^{(s,s/d)}_{p}\text{ when $s$ and }s/d\in
\Bbb N,\\
W^{(s,s/d)}_{p}&=B^{(s,s/d)}_{p}\text{ when $s$ and }s/d\in
\Bbb R_+\setminus\Bbb N.
\endaligned\tag A.11 
$$
(Note that the $W^{(s,s/d)}_p$ spaces do not always
interpolate well. Moreover, when e.g.\ $s\in\Bbb N$, $s/d\notin \Bbb N$,
one has a space where the $x$-regularity is $H^s_p$-type, the
$t$-regularity is $B^{s/d}_p$-type; for $p\ne 2$ it is not one of the
spaces in (A.10), cf.\ also (1.6).)

Observe that differential operators have the effect, for any
$s\in\Bbb R$,
$$\aligned
D_x^\alpha D_t^j\:H^{(s,s/d)}_p(\Rn\!\times\!\Bbb R)&\to
H^{(s-m,(s-m)/d)}_p(\Rn\!\times\!\Bbb R)\;\text{and } \\
D_x^\alpha D_t^j\:B^{(s,s/d)}_p(\Rn\!\times\!\Bbb R)&\to
B^{(s-m,(s-m)/d)}_p(\Rn\!\times\!\Bbb R),\;\text{for }|\alpha
|+dj\le m.
\endaligned\tag A.12
$$ 
This follows for the Bessel-potential spaces by Lizorkin's criterion
(cf.\ (2.11))
applied to $\xi ^\alpha \tau ^j\{\xi ,\tau \} ^{-m}$; it is seen
for the Besov spaces e.g\. from the definition by difference norms, or
by interpolation.

The spaces are defined over open subsets of ${\Bbb R}^{n+1}$ by
restriction; here the cylindrical subsets $\Sigma =\Omega \times I$ with
$\Omega $ open $\subset {\Bbb R}^n$ and $I$ an open interval of ${\Bbb
R}$ are particularly interesting. We use the notation  
$$
\aligned
&\ol H^{(s,s/d)}_p(\Omega \times I)=r_{\Omega \times I}H^{(s,s/d)}_p(\Rn\!\times\!\Bbb
R),
\quad
\ol B^{(s,s/d)}_p(\Omega \times I)=r_{\Omega \times I}B^{(s,s/d)}_p(\Rn\!\times\!\Bbb
R),
\\
&H^{(s,s/d )}_{p,\operatorname{loc}}(\Omega \times I )=\{u\in \D'(\Omega
 \times I )\mid \psi u\in  \ol H_p^{(s,s/d )}(\Omega \times I )\text{ for any
 }\psi \in C_0^\infty (\Omega \times I)\}, 
\endaligned\tag A.13
$$
and similar spaces $B^{(s,s/d )}_{p,\operatorname{loc}}(\Omega \times
   I ) $, $\ol W^{(s,s/d)}_p(\Omega \times I)$,  $W^{(s,s/d )}_{p,\operatorname{loc}}(\Omega \times I )$. Much of the above information, e.g.\ (A.7), (A.8), (A.10), carries
   over to the scales in the first line.

Let us also define anisotropic H\"older spaces. For  $k\in{\Bbb N}_0$, $0<\sigma <1$, the usual H\"older
   space of order $s=k+\sigma $ over
   $\Omega \subset{\Bbb R}^n$, in our notation $\ol C^s(\Omega )$, is
   provided with the norm
$$
\|u\|_{\ol C^s(\Omega )}={\sum}_{|\alpha |\le k}\|D^\alpha u\|_{L_\infty }+{\sum}_{|\alpha |= k}\operatorname{sup}_{x,x'\in\Omega}\frac{|D^\alpha u(x)-D^ au(x')|}{|x-x'|^\sigma }.
$$
Integer values of $s$ will be included in the scale by defining
$\ol C^k(\Omega )$  for $k\in{\Bbb N}_0$ 
 as the space of bounded
continuous functions with bounded continuous derivatives up to order
$k$ on $\comega$. (Then we have $\ol C^{s}(\Omega )\subset \ol
C^{s_1}(\Omega )$ whenever $s>s_1\ge 0$. There is a more refined choice of slightly
larger spaces for $k\in{\Bbb N}_0$, the H\"older-Zygmund spaces, that
fits better with interpolation theory, but which we shall not need
here, since we only show results ``with a loss of $\varepsilon $''.)

We then define the anisotropic H\"older space over $\Omega \times I$ ($I\subset {\Bbb
R}$) for $s>0$ %is then defined for noninteger positive exponents 
by:
$$
\ol C^{(s,s/d)}(\Omega \times I)=L_\infty (I,\ol C^s(\Omega))\cap
 \ol C^{s/d}( I, L_\infty (\Omega )), \text{ when }s>0. %\text{ $s$ and
%$s/d \notin{\Bbb N}$}.
\tag A.14 
$$
(Note that $ \ol C^{s/d} ( I, L_\infty (\Omega ))=  L_\infty (\Omega
   , \ol C^{s/d}( I))$.) In view of the well-known embedding properties for isotropic spaces,
   we have with $\varepsilon >0$:
$$\aligned
&\ol C^{(s,s/d )}(\Omega \times I)\subset\ol H^{(s-\varepsilon
   ,(s-\varepsilon )/d)}_p(\Omega \times I)\text{ and }\ol B^{(s-\varepsilon ,(s-\varepsilon )/d)}_p(\Omega \times
   I),             \text{ when }s>0,\\
&\ol H^{(s,s/d)}_p(\Omega \times I)\text{ and }\ol B^{(s,s/d)}_p(\Omega \times
   I)\subset \ol C^{(s-n/p-\varepsilon ,(s-n/p-\varepsilon )/d)}(\Omega \times
   I),\text{ if }s-n/p-\varepsilon >0;
\endaligned\tag A.15$$
%when the indices occurring in the H\"older
%spaces are noninteger; 
in the first inclusion we assume $\Omega \times
I$ to be bounded. Both inclusions are shown by comparing the spaces  via (A.10) and
(A.14). The first inclusion  follows immediately from the
isotropic case. For the second inclusion we can use that for small positive $\varepsilon _1<s$:
$$
\aligned
\ol H^{(s,s/d)}_p(\Omega \times I)&\subset \ol H_p^{\varepsilon _1
/d}(I,\ol H_p^{s-\varepsilon _1}(\Omega ))\cap  \ol
H_p^{(s-\varepsilon _1)/d}(I,\ol H_p^{\varepsilon _1
}(\Omega ))
\\
&\subset  L_\infty (I,\ol H_p^{s-\varepsilon _1}(\Omega ))\cap  \ol H_p^{(s-\varepsilon _1)/d}(I,L_\infty
(\Omega )),
\endaligned
$$
in order to relate to (A.14) (and similarly for $B$-spaces); this is allowed since $\ang\tau
^{\varepsilon _1/d}\ang\xi ^{s-\varepsilon _1}$ and $ \ang\tau
^{(s-\varepsilon _1)/d}\ang\xi ^{\varepsilon _1}$ are $\le c(\ang\xi
^{2d}+\tau ^2)^{s/(2d)}$.

There is also the local version:
$$
 C^{(s,s/d )}_{\operatorname{loc}}(\Omega \times I)=\{u\in \D'(\Omega
 \times I)\mid \psi u\in  \ol C^{(s,s/d )}(\Omega \times I)\text{ for any
 }\psi \in C_0^\infty (\Omega \times I)\}. \tag A.16
$$


\Refs

\widestnumber\key{[LPPS15]}
\ref
\no  [B64] 
\by O. Besov
\paper Investigation of a family of function spaces in connection
with theorems of imbedding and extension
\jour Trudy Mat. Inst. Steklov
\vol 60
\yr 1961
\pages 42--81
\moreref 
\jour = AMS Transl.
\vol 40
\yr1964
\pages85--126
\endref

\ref
\no  [BIN78] 
\by O. Besov, Il'in and Nikolski\u\i{}
\book Integral representations of functions and imbedding theorems
I--II
\publ Winston, Wiley
\publaddr Washington D. C.
\yr 1978 
\endref

\ref\no[BWZ17]\by U. Biccari, M. Warma and E. Zuazua \paper Local
regularity for fractional heat equations \finalinfo arXiv: 1704.07562
\endref

\ref\no[BBC03]\paper     Censored stable processes
 \by   K. Bogdan, K. Burdzy and Z.-Q. Chen \vol 127\pages 89--152
\jour Prob. Theory Related Fields\yr 2003 
\endref

\ref 
\no  [BM71] 
\by L. Boutet de Monvel
\paper Boundary problems for pseudo-differential operators
\jour Acta Math.
\vol 126
\yr1971
\pages11--51
\endref

\ref\no[C63] \by A. Calder\'on \paper Intermediate spaces and
interpolation \jour Studia Math., Seria specjalna\vol 1 \yr 1963
\pages 31--34
\endref

\ref\no[CD14]
\by H. Chang-Lara and G. Davila \paper Regularity for solutions of non
local parabolic equations\jour Calc. Var. Part. Diff.
Equations \vol 49 \yr2014 \pages 139-–172\endref

\ref\no[D89] \by E. B. Davies\book Heat kernels and spectral
theory. Cambridge Tracts in Mathematics, 92\publ Cambridge University
Press \publaddr Cambridge \yr 1989 
\endref

\ref\no[FK13]
\by M. Felsinger and M. Kassmann \paper Local regularity for parabolic
nonlocal operators\jour Comm. Part. Diff. Equations \vol 38
\yr2013 
\pages 1539-–1573\endref

\ref\no[FR17] 
\paper Regularity theory for general stable operators: parabolic equations
\by
X. Fernandez-Real and X. Ros-Oton\jour
J. Funct. Anal. \vol 272 \yr2017 \pages 4165--4221
\endref

\ref \no   [F69] 
\by A. Friedman 
\book Partial differential equations
\yr1969 
\publ Holt, Rinehart and Winston
\publaddr New York
\endref


\ref\no[FOT94]\by M. Fukushima, Y. Oshima and M. Takeda \book
Dirichlet forms and symmetric Markov processes. De Gruyter Studies in
Mathematics, 19\publ Walter de Gruyter \& Co.\publaddr Berlin \yr 1994
\endref

\ref \no   [G66] 
\by P. Grisvard 
\paper Commutativit\'e de deux foncteurs d'interpolation et
applications
\yr1966 
\vol 45 
\jour J. Math. Pures et Appl.
\pages 143--290
\endref


\ref \no   [G69] 
\by P. Grisvard 
\paper Equations diff\'erentielles abstraites
\yr1969 
\vol 2 (S\'erie 4)
\jour Ann. Ecole Norm\. Sup.
\pages 311--395
\endref

\ref\no[G90] \by G. Grubb \paper Pseudo-differential boundary problems
in $L_p$-spaces \jour Comm. Part. Diff. Equations \vol 13 \yr 1990 \pages
289--340
\endref

 \ref\no[G86]\by 
{G.~Grubb}\book Functional calculus of pseudodifferential
     boundary problems.
 Pro\-gress in Math.\ vol.\ 65 \publ  Birkh\"auser
\publaddr  Boston \yr 1986
\endref

\ref\no[G95]
\by G. Grubb\paper Parameter-elliptic and parabolic pseudodifferential
boundary problems in global Lp Sobolev spaces \jour Math. Z. \vol 218
\yr 1995 \pages 43-90
\endref

 \ref\no[G96]\by 
{G.~Grubb}\book Functional calculus of pseudodifferential
     boundary problems.
 Pro\-gress in Math.\ vol.\ 65, Second Edition \publ  Birkh\"auser
\publaddr  Boston \yr 1996
\endref


\ref\no[G09]\by G. Grubb\book Distributions and operators. Graduate
Texts in Mathematics, 252 \publ Springer \publaddr New York\yr 2009
 \endref


\ref\no[G14] \by G. Grubb \paper 
Local and nonlocal boundary conditions for $\mu $-transmission
and fractional elliptic pseudodifferential operators 
\jour Analysis and P.D.E. \vol 7 \yr 2014 \pages 1649--1682
\endref


\ref\no[G15] \by G. Grubb \paper Fractional Laplacians on domains, 
a development of H\"o{}rmander's theory of $\mu$-transmission
pseudodifferential operators \jour Adv. Math. \vol 268 \yr 2015
\pages 478--528
\endref

\ref\no[G17]\by G. Grubb \paper  \paper Green's formula and a Dirichlet-to-Neumann operator for
fractional-order pseudodifferential operators \finalinfo 
arXiv:1611.03024, to appear
\endref

\ref
\no [GK93] 
\by G\. Grubb and N\. J\. Kokholm
\paper A global calculus of parameter-dependent pseudodifferential
boundary problems in 
$L_p$ Sobolev spa\-ces 
\jour Acta Math.
\vol 171
\yr 1993
\pages 165--229
\endref


\ref
\no  [GS90] 
\by G. Grubb and V. A. Solonnikov
\paper Solution of parabolic pseu\-do-\-dif\-fe\-ren\-ti\-al in\-i\-ti\-al-boun\-d\-a\-ry
value problems
\jour J. Diff. Equ. 
\vol 87
\yr1990
\pages 256--304
\endref

\ref\no[H65]\by L. H\"o{}rmander\book Ch.\ II, Boundary problems for
``classical'' pseudo-differential operators  \finalinfo unpublished
lecture notes at Inst. Adv. Studies, Princeton 1965; available at
 http://www.math.ku.dk \linebreak /$\sim$grubb/LH65.pdf
\endref


\ref\no[H66]\by L. H\"ormander \paper Pseudo-differential operators and
non-elliptic boundary problems\jour Ann. of Math. \vol 83 \yr1966
\pages 129-209 
\endref

\ref\no[H85]\by L. H\"o{}rmander\book The analysis of linear partial
differential operators, III \publ Springer Verlag\publaddr Berlin, New
York\yr 1985
 \endref

\ref\no[JX15]\by 
T. Jin and J. Xiong \paper Schauder estimates for solutions of linear
parabolic integro-differential equations \jour Discrete
Contin. Dyn. Syst. \vol 35 \yr2015 \pages 5977-–5998\endref

\ref\no[KN65]\by J.J. Kohn and L. Nirenberg\paper An algebra of
pseudo-differential operators \jour Comm. Pure Appl. Math.\vol 18 \yr
1965\pages  269-305
\endref

\ref \no   [LSU68] 
\by O. A. Ladyzhenskaya, V. A. Solonnikov and N. N. Uraltzeva
\book Linear and Quasilinear Equations of Para\-bo\-lic Type
\yr1968 
\publ Amer. Math. Soc.
\publaddr Providence
\endref


\ref\no[L87]\by D. Lamberton\paper \'E{}quations d'\'evolution
lin\'eaires associ\'ees a des semi-groupes de contractions dans les
espaces Lp \jour J. Funct. Anal. \vol 72 \yr 1987 \pages 252-262
\endref

\ref\no[LPPS15]\by T. Leonori, I. Peral, A. Primo and F. Soria \paper
Basic estimates for solutions of a class of nonlocal elliptic and
parabolic equations
\jour
Discrete Contin. Dyn. Syst. \vol 35 \yr 2015 \pages 6031-6068
\endref



\ref \no   [LM68] 
\by J.-L. Lions and E. Magenes 
\book Probl\`emes aux limites non homog\`enes et applications. Vol. 1
et 2
\yr1968 
\publ Editions Dunod 
\publaddr Paris
\endref


\ref
\no  [L66] 
\by P. I. Lizorkin
\paper Nonisotropic Bessel potentials. Imbedding theorems for Sobolev
spaces  $L_p^{r_1,\dots,r_n}$, fractional derivatives.
\jour Dokl. Akad. Nauk SSSR
\vol 170
\yr 1966
\pages 508--511
\moreref 
\jour = Soviet Math. Doklady
\vol 7
\yr1966
\pages1222--1225
\endref


\ref
\no  [L67] 
\by P. I. Lizorkin
\paper Multipliers of Fourier integrals in the space $L_p$
\jour Trudy Mat. Inst Steklov
\vol 89
\yr 1967
\pages 247--267
\moreref 
\jour = Proc. Steklov Inst. Math.
\vol89
\yr1967
\pages269--290
\endref

\ref 
\no  [N75] 
\by S. M. Nikolski\u\i
\book Approximation of functions of several variables and imbedding
theorems
\publ Springer Verlag
\publaddr Berlin
\yr 1975
\endref


\ref\no[R16]\paper
Nonlocal elliptic equations in bounded domains: a survey
\by X. Ros-Oton
\jour Publ. Mat. \vol60 \yr2016 \pages 3--26\endref

\ref 
\no  [SC66] 
\by C. Sadosky and M. Cotlar
\paper On quasi-homogeneous Bessel potential operators
\jour Proc. Symp. Pure Math. 
\vol 10
\yr 1966
\pages 275--287
\endref


\ref\no[S67] \by M. Schechter \paper Complex interpolation \jour
Compositio Math. \vol 18 \yr 1967 \pages 117--147
\endref

\ref\no[ST76]\by H.-J. Schmeisser and H. Triebel\paper Anisotropic
spaces I. (Interpolation of abstract spaces and function spaces)
\jour Math. Nachr. \vol73 \yr 1976 \pages 107--123
\endref


\ref
\no  [ST87] 
\by H.-J. Schmeisser and H. Triebel
\book Fourier analysis and function spaces
\publ Wiley
\publaddr New York
\yr1987
\endref


\ref\no[S65]\by R. T. Seeley
\paper Integro-differential operators on vector bundles \jour
Trans. Amer. Math. Soc. \vol 117 \yr 1965\pages 167-204
\endref

\ref
\no   [So65] 
\by V. A. Solonnikov
\paper On boundary-value problems for linear parabolic systems of
differential equations of general form 
\jour Trudy Mat. Inst. Steklov
\vol 83
\yr 1965
\pages 1--162
\moreref
\jour = Proc. Steklov Inst. Math.
\vol 83
\yr 1965
\pages 1--162
\endref

\ref
\no  [S87] 
\by V. A. Solonnikov
\paper On the transient motion of an isolated volume of viscous
incompressible fluid
\jour Izv. Akad. Nauk SSSR
\vol51
\yr 1987
\pages1065--1087
\moreref
\jour = Math. USSR Izvestiya
\vol31
\yr1988
\pages 381--405
\endref

\ref\key[T81]\by M.~E. Taylor\book
 Pseudodifferential operators \publ
Princeton University Press \publaddr Princeton, NJ \yr1981
\endref


\ref
\no  [T78] 
\by H. Triebel
\book Interpolation theory, function spaces, differential operators
\publ North-Hol\-land Publ. Co.
\publaddr Amsterdam, New York
\yr1978
\endref



\ref\no[Y84]
\by M. Yamazaki\paper The Lp-boundedness of pseudodifferential
operators satisfying estimates of parabolic type and product type\jour
Proc. Japan Acad. Ser. A Math. Sci. \vol 60 \yr 1984 \pages 279-282
\endref

\ref\no[Y86] \by
M. Yamazaki\paper A quasihomogeneous version of paradifferential
operators. I. Boundedness on spaces of Besov type \jour
J. Fac. Sci. Univ. Tokyo Sect. IA Math. \vol 33 \yr 1986 \pages 131-174 
\endref


\endRefs
\enddocument